\documentclass[11pt]{article}
\usepackage{tipa}
\usepackage[all,import]{xy}
\usepackage{graphicx,color}
\usepackage{amsmath}
\usepackage{amsbsy}
\usepackage{amssymb}
\usepackage{cases}
\usepackage{enumerate}
\usepackage{bm}
\usepackage{amsthm}
\usepackage{multirow}
\usepackage{subfigure}
\usepackage{colortbl}
\definecolor{mygray}{gray}{.9}

\usepackage[numbers, sort&compress]{natbib}
\makeatletter 
\renewcommand\@biblabel[1]{#1.} 
\makeatother %

\begin{document}
\theoremstyle{definition}
\newtheorem{assumption}{Assumption}
\newtheorem{theorem}{Theorem}
\newtheorem{lemma}{Lemma}
\newtheorem{example}{Example}
\newtheorem{definition}{Definition}
\newtheorem{corollary}{Corollary}

\def\letas{\mathrel{\mathop{=}\limits^{\triangle}}}
\def\ind{\begin{picture}(9,8)
         \put(0,0){\line(1,0){9}}
         \put(3,0){\line(0,1){8}}
         \put(6,0){\line(0,1){8}}
         \end{picture}
        }
\def\nind{\begin{picture}(9,8)
         \put(0,0){\line(1,0){9}}
         \put(3,0){\line(0,1){8}}
         \put(6,0){\line(0,1){8}}
         \put(1,0){{\it /}}
         \end{picture}
    }

\def\var{\text{var}}
\def\cov{\text{cov}}
\def\sumn{\sum\limits_{i=1}^n}
\def\summ{\sum\limits_{j=1}^m}
\def\convergeas{\stackrel{a.s.}{\longrightarrow}}
\def\converged{\stackrel{d}{\longrightarrow}}
\def\iidsim{\stackrel{i.i.d.}{\sim}}
\def\indsim{\stackrel{ind}{\sim}}
\def\asim{\stackrel{a}{\sim}}

\setlength{\baselineskip}{1.5\baselineskip}

\begin{flushleft}
{\Large \bf Identifiability of Subgroup Causal Effects in Randomized Experiments with Nonignorable Missing Covariates}\\

Peng Ding$^1$ and Zhi Geng$^{2}$  \\
1 Department of Statistics, Harvard University, Science Center, One Oxford Street, Cambridge, MA 02138, U.S.A. (E-mail: \texttt{pengding@fas.harvard.edu})\\
2 Center for Statistical Science, 
School of Mathematical Sciences, 
Peking University, Beijing
100871, China (E-mail: \texttt{zhigeng@pku.edu.cn})\\

\end{flushleft}

\begin{quotation}
\noindent

\noindent {\bf Summary}. \quad
Although randomized experiments are widely regarded as the gold standard for estimating causal effects, missing data of the pretreatment covariates makes it challenging to estimate the subgroup causal effects.
When the missing data mechanism of the covariates is nonignorable, the parameters of interest are generally not pointly identifiable, and we can only get bounds for the parameters of interest, which may be too wide for practical use.
In some real cases, we have prior knowledge that some restrictions may be plausible.
We show the identifiability of the causal effects and joint distributions for four interpretable missing data mechanisms, and evaluate the performance of the statistical inference via simulation studies.
One application of our methods to a real data set from a randomized clinical trial shows that one of the nonignorable missing data mechanisms fits better than the ignorable missing data mechanism, and the results conform to the study's original expert opinions.
We also illustrate the potential applications of our methods to observational studies using a data set from a job-training program.

\end{quotation}

\section{Introduction}
Randomized experiments are widely regarded as the gold standard for estimating causal effects.  However, a key problem in the analysis is missing data due to ethical or logistical reasons, either pretreatment covariates [1-6]
or outcomes [3, 7-9].
In general, without making untestable assumptions, we can only obtain large sample bounds for parameters of interest, rather than point estimates [3, 5, 6].
In practice, however, these bounds may be too wide for practical use.
If the missing data mechanism is ignorable,
the likelihood and Bayesian inference based on the observed data are both valid.
In many cases, however, the missing data mechanism is not ignorable; that is, the missing data process depends on some possibly missing variables.
Sensitivity analysis was used in nonignorable missing covariates problem [5, 10],
however, sensitivity analysis could not provide point identification.
There are some works on identification analysis when the missingness of the outcomes were nonignorable [7, 9, 11].
Previous works discussed the identifiability of causal effects when a key covariate is truncated due to death [1, 2, 6].
And nonignorable missing covariates problems were also discussed in survival analysis [4]
and regression models [12].

In this paper, we discuss nonignorable missing data mechanisms of a key covariate.
Although the average treatment effect for the whole population is still identifiable by randomization,
a central research question concerns the subgroup causal effects defined by the possibly missing covariate.
Because of ethical and logistical problems in some randomized experiments, there may be strong evidence that the missing data may depend on the missing covariates.
However, it is quite challenging to make inference about the nonignorable missing data problem, because the joint distribution is not identifiable without assumptions or restrictions.
Bounds of the subgroup treatment effects were obtained with and without some expert opinions using numerical optimization methods [5]. 
We take a slightly different approach from the bound analysis, and instead focus exclusively on making assumptions about the missing data mechanism itself.
Based on plausible assumptions and restrictions, we theoretically demonstrate the identifiability of subgroup causal effects, and obtain explicit forms for the bounds of the parameters of interest. 
The simulation study shows promising results about the finite sample performance of our methods.
We first apply the proposed methods to a real randomized clinical trial. Interestingly,
we conclude that one of the nonignorable missing data mechanisms fits the data better than the ignorable missing data mechanism, and the corresponding point estimates satisfy
all the experts' ``scientific assumptions'' proposed in the original analysis [5].
We also analyze a job-training program data as an application of our methods to an observational study.

The rest of the paper is organized as follows.
In Section \ref{sec::notation}, we introduce the notation and the main assumptions used throughout the paper.
In Section \ref{sec::missing}, we discuss possible missing data mechanisms for a pretreatment covariate of interest.
In Section \ref{sec::identifiability}, we establish the identifiability of these models.
In Section \ref{sec::compute}, we describe the computational details and related hypothesis testing problems. 
In Section \ref{sec::simulation}, we validate these findings using a simulation study.
In Section \ref{sec::real}, we turn to a randomized clinical trial, and address both model selection and the assumptions from the original study.
In Section \ref{sec::job-train}, we analyze a data set from a job-training program.
We conclude with a discussion, and present the details of the computations, simulations and proofs in the web Appendix.

\section{Notation and Assumptions}\label{sec::notation}
We are interested in a randomized experiment with $N$ subjects.
Suppose that $T$ is a binary treatment variable.
For a subject $i$, let $T_i$ denote the treatment assignment, with
$T_i=1$ if subject $i$ receives the treatment,
and $T_i=0$ if subject $i$ receives the control.
Let $X_i$ denote the pretreatment covariate with $J$ categories,
which may be missing.
Let $Y_i$ denote the observed outcome with $K$ categories.
Furthurmore,
let $M_i$ be the missing data indicator for $X_i$,
with $M_i = 1$ if $X_i$ is missing, and $M_i = 0$ if $X_i$ is observed.

We define causal effects via the potential outcomes model.
Suppose that missing of $X$ happens after the treatment assignment.
Thus both the outcome and the missing of $X$ may be affected by treatment assignment $T=t$.
Let $Y_i(t)$ and $M_i(t)$ denote the potential outcome variable and the potential missing data indicator for $X_i$ respectively, if subject $i$ were assigned to treatment $t$.
These variables are potential outcomes because only one of the pairs $\{  Y_i(1), M_i(1)\} $ and $\{ Y_i(0), M_i(0)\}$ can be observed.
Let $Y_i$ and $M_i$ denote their observations.
The potential outcomes are well-defined under the following fundamental and widely used assumptions in causal inference.

\begin{assumption} \label{assum::sutva}
(Stable unit treatment value assumption, SUTVA)
There is no interference between units, which means that the potential outcomes
of one individual do not depend on the treatment status of other
individuals [13],
and there is only one version of the potential outcome of a certain treatment [14].
\end{assumption}

The SUTVA assumption implies that the observed outcomes are deterministic functions of the potential outcomes and the treatment assignments, i.e., $Y_i = Y_i(T_i)  = T_iY_i(1) + (1 - T_i)Y_i(0)$ and $M_i = M_i (T_i) = T_i M_i(1) + (1 - T_i) M_i(0)$.
Throughout our paper, we assume that $\{(T_i, Y_i(1), Y_i(0), M_i(1), M_i(0), X_i): i =1, ..., N\}$ are random draws from a superpopulation, and therefore the observed data $\{ (T_i, Y_i, X_i, M_i ): i=1,...,N\}$ are also randomly drawn from the superpopulation, where $X_i$ is missing when $M_i=1.$

Let $A\ind B | C$ denote that variables $A$ and $B$ are conditionally independent given variable $C$.
The following assumption is a generalization of the ignorable treatment assignment assumption [15].

\begin{assumption}  \label{assum::random}
(Latent ignorable treatment assignment)
$T\ind \{Y(1), Y(0), M(1), M(0) \}|X$.
\end{assumption}

When there is no missing data in $X$, Assumption \ref{assum::random} is called the ``ignorable assumption'' of the treatment assignment mechanism, which is crucial for causal inference in observational studies.
When $X$ has some missing values, the assumption is no longer the original ``ignorable treatment assignment'' mechanism, and the difficulty for point identification of the causal effects arises.
As pointed out by a reviewer, it is a ``latent ignorable treatment assignment'', since the observed data does not contain all the values of $X$. Previous researchers [11] 
first used this term for nonignorable missing data of the outcomes, and also used it for treatment assignment mechanism with some key covariates missing [16].

In completely randomized experiments, the following stronger assumption holds by the design of experiments.

\begin{assumption}\label{assum:strong_random}
 (Complete randomization)  
 $T\ind \{Y(1), Y(0), M(1), M(0), X  \}$.
 \end{assumption}

The independence assumption above means that $P\{ T\mid Y(1), Y(0), M(1), M(0), X \}  = P(T)$, under the assumption of existence of a joint distribution of $\{ T, Y(1), Y(0), M(1), M(0), X\}$. The joint independence assumption above implies Assumption \ref{assum::random}.
Since all the potential outcomes and $X$ are ``pretreatment covariates'', the treatment assignment mechanism $T$ is independent of all of them in completely randomized experiments.
In one of our real applications in this paper, the data comes from a completely randomized experiment, and the stronger Assumption \ref{assum:strong_random} is satisfied automatically. However, the theory and methods discussed in this paper can be applied to more general problems under a weaker Assumption \ref{assum::random}, which may be more plausible in observational studies. We also use another example to illustrate the potential applications of our method in observational studies.

The goal of this paper is to use the observed data to make inference about the following measure of causal effects:
$$
CE_x =  \mathcal{D}\left[ E\{ Y(1)  \mid X=x \}, E \{ Y(0) \mid X=x\} \right],
$$
where
$\mathcal{D}[p_1, p_0]$ is a function with the following properties:
\begin{enumerate}[1.~]
\item
$\partial\mathcal{D}/\partial p_1 >0$, $\partial\mathcal{D}/\partial p_0 <0$, and

\item
$\mathcal{D}[p_1, p_0]$ and $p_1 - p_0$ have the same sign,
$>0$, $<0$ or $=0$.
\end{enumerate}
We say that a treatment $T$ has a positive (negative or null) causal effect
on an outcome $Y$ in subgroup with $X=x$, if $CE_x$ is larger than (smaller or equal to) zero.
For example, for a binary $Y$,
let $p_1=P(Y(1) =  1 \mid X=x)$ and $p_0=P(Y(0) = 1\mid X=x)$.
Then $CE_x$ may be
the causal risk difference (CRD) $\mathcal{D}[p_1, p_0] = p_1 - p_0$,
the log of the causal risk ratio (CRR) $\log(p_1/p_0)$ or the log of the causal odds ratio (COR) $\log [ p_1(1-p_0)/\{p_0(1-p_1)\} ]$.
The total causal effect $CE_+ =  \mathcal{D}\left[   E\{ Y(1) \} , E\{ Y(0) \} \right]$ is also of interest in practice.

Under Assumption \ref{assum:strong_random}, $CE_+$ is identifiable because $P\{ Y(t)=y \} = P(Y=y\mid T=t)$ for $t=0$ and $1.$
Under Assumption \ref{assum::random}, we have $P\{ Y(t)=y \mid X=x \} = P\{ Y(t) =y \mid T=t, X=x\} = P(Y=y \mid T=t, X=x)$.
The causal effects can be expressed as functions of the joint distribution of $(T,X,Y)$, and
$CE_x$ is identifiable if $P(T=t,X=x,Y=y)$ is identifiable.

Define $p_{txym} = P(T=t, X=x, Y=y, M=m)$ and $p_{t+y1} = P(T=t, Y=y, M=1)$,
where ``+'' in a subscript denotes the distribution marginalized over the corresponding variable.
Analogously, let $N_{txym}$ denote the observed frequency in
the cell $(t, x, y, m)$ of the contingency table,
and let $N_{t+ym}$ denote the marginal frequency of the contingency table over
the corresponding variable $X$.
We can directly identify $p_{txy0}$ and $p_{t+y1}$ by the observed frequencies, $N_{txy0}/N$ and $N_{t+y1}/N$, respectively.
However, we cannot identify $p_{txy1}$ and thus $CE_x$ without any further assumptions.
Throughout the paper, we need the following condition that the
data is not missing with probability one.

\begin{assumption}
\label{assum::missing}
$P(M=0\mid T=t, X=x, Y=y)>0$ for all $t,x,y$.
\end{assumption}


\section{Missing Data Mechanisms}\label{sec::missing}
Before discussing the missing data mechanisms, we first review the some definitions about missing data [17].
Let $\bm{D}_{com} $ be the generic notation for the complete data, 
$\bm{D}_{obs}$ for the observed data, $\bm{D}_{mis}$ for the missing data, and $\bm{D}_{com} = ( \bm{D}_{obs}, \bm{D}_{mis}  )$. Let $\bm{M}$ be the indicator matrix of the missing data corresponding to $\bm{D}_{com} $. The model for the complete data is $f(\bm{D}_{com}\mid \theta)$, and the model for $\bm{M}$ given $\bm{D}_{com}$ is $f(\bm{M}\mid \bm{D}_{com}, \psi)$. 

\begin{definition}
The missing data mechanism is called missing at random (MAR), if $\bm{M}$ only depends on $\bm{D}_{obs}$, i.e., $f(\bm{M}\mid \bm{D}_{com}, \psi) = f(\bm{M}\mid \bm{D}_{obs}, \psi)$. Otherwise, if $\bm{M}$ depends on $\bm{D}_{mis}$, the missing data mechanism is called missing not at random (MNAR).
\end{definition}

\begin{definition}
We call the parameters $\theta$ and $\psi$ distinct, if the parameter space of $(\theta, \psi)$ is the product of the parameter space of $\theta$ and the parameter space of $\psi$.
\end{definition}

\begin{definition}
The missing data mechanism is called ignorable, if it is MAR and the parameters $\theta$ and $\psi$ are distinct.
\end{definition}

The likelihood for $(\theta, \psi)$ is proportional to
$$
f(\bm{D}_{obs}, \bm{M}\mid \theta, \psi)  = \int f(\bm{D}_{obs}, \bm{D}_{mis} \mid \theta) f( \bm{M}\mid \bm{D}_{obs}, \bm{D}_{mis}, \psi ) d\bm{D}_{mis}.
$$
Under the ignorable missing data mechanism, it reduces to
$$
f(\bm{D}_{obs}, \bm{M}\mid \theta, \psi) = f(\bm{M}\mid \bm{D}_{obs}, \psi) f( \bm{D}_{obs}\mid \theta ).
$$
In this case, the inference for $\theta$ can be based only on the observed data likelihood $f( \bm{D}_{obs}\mid \theta )$, and the missing data mechanism can be ``ignored''.
For the missing data mechanisms discussed in this paper, the parameters are distinct, and therefore ``nonignorable'' is equivalent to NMAR.

In randomized experiments with $X$ missing, missing data mechanisms influence the identifiability and estimation of subgroup causal effects $CE_x$.
We consider five missing data mechanisms.
The first one is ignorable where the missing $X$ depends
only on observed variables $(Y,T)$.
The others are nonignorable.
For mechanisms \ref{miss::2} and \ref{miss::3},
we assume that the missingness of $X$ depends on both $X$
and another one of $Y$ and $T$.
For the last two mechanisms 4 and 5, we assume that they
depend on all three variables $(X,Y,T)$.
Under Assumption \ref{assum::random},
the five missing mechanisms to be discussed
in the next section can be described equivalently
in terms of both potential and observed outcomes as follows:
\begin{enumerate}[$M_1.~$]
\item\label{miss::1}
$M$ may depend on $(Y, T)$ but is independent of $X$ conditional on $(Y, T)$, i.e.,
\begin{eqnarray*}
P\{ M(t) = 1\mid X=x, Y(t)=y\}  &=& P\{ M(t)=1 \mid T=t, Y(t) =y \} , \text{ or }\\
P(M =1\mid T=t, X=x, Y=y)  &=& P(M=1\mid T=t, Y=y);
\end{eqnarray*}
further, if the randomization Assumption \ref{assum:strong_random} holds,
the missing mechanism is equivalent to
$P\{  M(t) = 1\mid X=x, Y(t)=y\}  = P\{ M(t)=1 \mid Y(t) =y \}$;

\item\label{miss::2}
$M$ may depend on $(X, T)$ but is independent of $Y$ conditional on $(X, T)$, i.e.,
\begin{eqnarray*}
P\{  M(t) = 1\mid X=x, Y(t)=y\}  &=& P\{  M(t)=1 \mid X=x\}, \text{ or } \\
P(M =1\mid T=t, X=x, Y=y)  &=& P(M=1\mid T=t,  X=x );
\end{eqnarray*}

\item\label{miss::3}
$M$ may depend on $(X, Y)$ but is independent of $T$ conditional on $(X, Y)$, i.e.,
\begin{eqnarray*}
P\{  M(t) = 1\mid X=x, Y(t)=y\}&=& P\{ M(t') = 1\mid X=x, Y(t')=y\} \text{ for $t\neq t'$, or }\\
P(M =1\mid T=t, X=x, Y=y)  &=& P(M=1\mid X=x, Y=y);
\end{eqnarray*}

\item\label{miss::4}
$M$ may depend on $(X, Y, T)$ via a Logistic model, i.e.,
\begin{eqnarray*}
\text{logit} [ P\{ M(t)=1\mid X=x,Y(t)=y \} ]
&=& \beta_0 + \beta_T t + \beta_X x+ \beta_Y y, \mbox{ or } \\
\text{logit} \{  P(M=1\mid T=t, X=x,Y=y) \}
&=& \beta_0 + \beta_T t + \beta_X x+ \beta_Y y,
\end{eqnarray*}
where $\text{logit}\{a\} = \log\{ a/(1-a)\}$;

\item\label{miss::5}
$M$ may depend on $(X,Y,T)$ and does not have any restrictions.

\end{enumerate}

\section{Identifiability of Causal Effects}\label{sec::identifiability}

In this section, we discuss the identifiability of causal effects and the joint distribution of $(X, Y, T, M)$ for the missing mechanisms
presented in Section \ref{sec::missing}.
If the joint distribution of $(X, Y, T, M)$ is identifiable, the causal effects $CE_x$ and $CE_+$ are also identifiable under Assumptions 2 or 3.
For the first four missing mechanisms,
we shall show that causal effects are identifiable,
and we shall give conditions for identifiability of the joint distribution of $(X, Y, T, M)$.
For the last missing mechanism,
we shall give lower and upper bounds for causal effects.

\begin{theorem}\label{thm::miss1}
For missing mechanism \ref{miss::1},  under Assumptions \ref{assum::sutva} and \ref{assum::random},
the joint distribution of $(T,X,Y,M)$ is identifiable.
\end{theorem}

The missing mechanism \ref{miss::1} is ignorable,
and the joint distribution can be consistently estimated from the observed data.

\begin{theorem}\label{thm::miss2}
For missing mechanism \ref{miss::2},
under Assumptions \ref{assum::sutva} and \ref{assum::random}, 

\noindent (1) the causal effects $CE_x$ are identifiable;

\noindent (2)
the joint distribution of $(T,X,Y,M)$ is identifiable if
Rank$\left( \bm{\Theta}_t \right) = J$ for $t=0$ and 1,
where $\bm{\Theta}_t$ is a $J \times K$ matrix with
$p_{txy0}$ as the $(x,y)$ element; and

\noindent (3)
for binary $X$, the rank condition reduces to $X\nind Y|(T=t)$ for $t=0$ and $1$, which is equivalent to the testable condition $X\nind Y|(T=t,M=0)$ for $t=0$ and $1$.
\end{theorem}

From Theorem \ref{thm::miss2},
we can also see that
under Assumptions \ref{assum::sutva} and \ref{assum::random} and missing mechanism 2,
$CE_x$ is always identifiable,
but the joint distribution may not be identifiable,
because the number of parameters is larger
than the number of observed frequencies
if $J>K$.
The rank condition for identifying the joint distribution can be checked, because the rank of $\bm{\Theta}_t$ equals the rank of the matrix with  $P(X=x, Y=y\mid T=t, M=0)$ as the $(x,y)$ element, which can be identified by the observed data.
It is necessary for the rank condition that $Y$ has more categories than $X$ (i.e., $J\leq K$) and
that there exists a subset of $Y$'s categories,
$\Omega = \{y_1,...,y_J \}\subseteq \{1,...,K\} $, such that
$P(X=x\mid T=t, M=0, Y=y) \neq P(X=x\mid T=t, M=0, Y=y')$ for any $y \neq y' \in \Omega$.

\begin{theorem}\label{thm::miss3}
For the missing mechanism \ref{miss::3},

\noindent (1)
under Assumptions \ref{assum::sutva} and \ref{assum::random}, if $Y$ is binary,
the log of causal odds ratios
$\log( COR_x )$ are identifiable, but only the signs of other causal effects $CE_x$ are identifiable;

\noindent (2)
under Assumptions \ref{assum::sutva} and \ref{assum:strong_random},  if $Y$ is binary,
the causal effects $CE_x$ and $CE_+$ are identifiable; and

\noindent (3)
under Assumptions \ref{assum::sutva} and \ref{assum:strong_random} or
under Assumptions \ref{assum::sutva} and \ref{assum::random},
the joint distribution of $(T,X,Y,M)$ is identifiable
if $X$ is binary and $X \nind T| (Y=y)$ for $y=0, 1, \ldots, K-1$.
The conditional dependence $X  \nind T| (Y=y)$ is equivalent to the testable condition
$X \nind T| (Y=y,M=0)$.
\end{theorem}

Since we compare only two treatment groups (i.e., $T$ is binary),
the condition for identifying the joint distribution for missing mechanism 3 requires that $X$ is binary;
otherwise the number of observed frequencies
are smaller than the number of parameters,
and thus the joint distribution is not identifiable.

The missing mechanism with $M \ind (T,Y)|X$ is a special case of both the missing mechanisms \ref{miss::2} and \ref{miss::3}, for which we have the following corollary from Theorems \ref{thm::miss2} and \ref{thm::miss3}.

\begin{corollary}
\label{corollary::M-X}
For the missing mechanism $M \ind (T,Y)|X$, under Assumptions \ref{assum::sutva} and \ref{assum::random}, the causal effects $CE_x$ are identifiable;
furthermore the joint distribution of $(T,X,Y,M)$ is identifiable if
Rank$\begin{pmatrix} \bm{\Theta}_1\\ \bm{\Theta}_0 \end{pmatrix} = J$. When $X$ is binary, the rank condition is equivalent to $X \nind (T,Y)$ which is further equivalent to the testable condition $X \nind (T,Y)|(M=0)$.
\end{corollary}

\begin{theorem}\label{thm::miss4}
Assume that $X$ and $Y$ are binary, and Assumptions \ref{assum::sutva} and \ref{assum::random} hold.
For the missing mechanism \ref{miss::4},
the joint distribution of $(T,X,Y,M)$ is identifiable
if the value of $OR_{YT|(M=1)}$
is between those of
$OR_{YT|(X=1,M=0)}$ and $OR_{YT|(X=0,M=0)}$,
where $OR_{YT|(M=1)}$ and $OR_{YT|(X=x,M=0)}$ are the
odds ratios of $Y$ and $T$ conditional on $M=1$ and on $(X=x, M=0)$, respectively.
\end{theorem}

Because $OR_{YT|(M=1)}$ and $OR_{YT|(X=x,M=0)}$
are identifiable,
the condition can be checked by the observed data. In our application to a randomized clinical trial in Section 7, the condition is satisfied.

\begin{theorem}\label{thm::miss5}
Assume that $Y$ is binary, and Assumptions \ref{assum::sutva} and \ref{assum::random} hold.
For the missing mechanism \ref{miss::5},
the lower bound for $CE_x$ is
$$
\mathcal{D}\left[  \frac{  p_{1x10}    }{  p_{1x00} + p_{1x10} + p_{1+01}     }   ,  \frac{  p_{0x10} + p_{0+11}  }{    p_{0x00} + p_{0x10} + p_{0+11}    }  \right],
$$
which is attainable when $p_{1x01} = p_{1+01}$, $p_{0x11 } = p_{0+11}$, and $p_{1x11} = p_{0x01} = 0$.
The upper bound for $CE_x$ is
$$
\mathcal{D}\left[ \frac{   p_{1x10} + p_{1+11}    }{   p_{1x00} + p_{1x10} + p_{1+11}   }      , \frac{p_{0x10}  }{p_{0x00}    + p_{0x10}   + p_{0+01}    }   \right],
$$
which is attainable when $p_{1x01} = p_{1+01}$, $p_{0x01} = p_{0+01}$ and
$p_{1x11} = p_{0111} = 0$.
\end{theorem}

The bounds of $CE_x$ given in Theorem \ref{thm::miss5} can be estimated from observed data by replacing the cell probabilities with the cell counts, but they may cover zero
and thus we may not be able to determine the sign of $CE_x$.
All estimates of $CE_x$ obtained under missing mechanisms \ref{miss::1} to \ref{miss::4} should fall into the bounds given in Theorem \ref{thm::miss5},
and particularly this can be shown for the maximum likelihood estimators (MLEs) of $CE_x$.
Under Assumptions \ref{assum::sutva} and \ref{assum:strong_random}, the bounds should be narrower but have no explicit forms, optimization methods can be used to find the numerical solutions.

\section{Computational Details and Hypothesis Testing Problems}
\label{sec::compute}
\subsection{EM Algorithms and Gibbs Samplers}
In practice, we can use the Expectation-Maximization (EM) algorithm to find the MLEs and use the Gibbs Sampler to simulate the posterior distributions of the parameters.
In this subsection, we only describe the computational details for missing mechanism 1, and the web Appendix provides more details for missing mechanisms 2 to 4.
For simplicity, we only describe the algorithms for binary $X$ and binary $Y$, and the algorithms for categorical $X$ and $Y$ can be written similarly. 
Denote $P^{(j)}(T=t, X=x, Y=y, M=m) = P^{(j)}(X=x) P^{(j)}(T=t\mid X=x)  P^{(j)}(Y=y\mid T=t, X=x) P^{(j)}(M= m \mid T=t, X=x)$ as the joint distribution of $(T,X,Y,M)$ in the $j$-th iteration for either the EM algorithm or the Gibbs Sampler. Define
\begin{eqnarray*}
 p_{x|ty1}^{(j)} &=& P^{(j)}(X=x\mid T=t, Y=y, M=1) \\
 &=& \frac{ P^{(j)}(T=t, X=x, Y=y, M=1) }  
 {  \sum_{x'=0,1} P^{(j)}(T=t, X=x', Y=y, M=1)  }.
\end{eqnarray*}
The EM algorithm iterates between the following E-step and M-step:
\begin{itemize}
\item
E-step: The sufficient statistics are imputed as $N^{(j)}_{txy0} = N_{txy0}$ and $N^{(j)}_{txy1} = N_{t+y1} p_{x|ty1}^{(j)} $;

\item
M-step: 
The joint distribution is updated by 
$P^{(j+1)}(T=t, X=x, Y=y, M=m) = \frac{N^{(j)}_{+x++}}{N^{(j)}_{++++}}   \frac{N^{(j)}_{tx++}}{N^{(j)}_{+x++}} 
 \frac{N^{(j)}_{txy+}}{N^{(j)}_{tx++}}  \frac{N^{(j)}_{tx+m}}{N^{(j)}_{tx++}}$.

\end{itemize}

\noindent The Gibbs sampler iterates between the following Imputation-step and Posterior-step:
\begin{itemize}
\item
Imputation-step: 
We let $N^{(j)}_{txy0} = N_{txy0}$ and draw $N^{(j)}_{txy1} \sim \text{Binomial}(N_{t+y1}, p_{x|ty1}^{(j)}) $;

\item
Posterior-step:
We draw 
$ P^{(j+1)}(X=1) \sim \text{Beta} (\alpha_X + N^{(j)}_{+1++}, \beta_X +N^{(j)}_{+0++} )$,
$P^{(j+1)}(T=1\mid X=x)\sim \text{Beta}( \alpha_T^x + N^{(j)}_{1x++}, \beta_{T}^x + N^{(j)}_{0x++}  )$,
$ P^{(j+1)}(Y=1\mid T=t, X=x) \sim \text{Beta} (\alpha_{Y}^{tx} + N^{(j)}_{tx1+} , \beta_{Y}^{tx} + N^{(j)}_{tx0+} )$,
$  P^{(j+1)}(M=1\mid T=t, X=x) \sim \text{Beta} (\alpha_{M}^{tx} + N^{(j)}_{tx+1} , \beta_{M}^{tx} + N^{(j)}_{tx+0} )$,
where $\alpha_X, \alpha_T^x, \beta_T^x,  \beta_X, \alpha_Y^{tx}, \beta_{Y}^{tx}, \alpha_M^{tx}$, and $ \beta_M^{tx}$ are parameters for the Beta priors of the probability parameters.
\end{itemize}

In our simulation studies and applications, we use the conventional noninformative Beta$(1/2, 1/2)$ prior for the probability parameters, the results of which are similar to the results from another commonly-used Uniform$(0,1)$ prior, when the sample sizes are relatively large.

\subsection{Testing Goodness-of-Fit, Interaction and Effect Modification}
\label{sec::test}

Under Assumption \ref{assum::random}, none of missing mechanisms \ref{miss::1} to \ref{miss::4} are testable from the observed data, since the numbers of parameters are equal to the numbers of observed frequencies.
Under Assumption \ref{assum:strong_random}, we have an additional constraint that $T\ind X$ by complete randomization. 
For example, in our first application with a binary outcome $Y$ and a binary covariate $X$, we can perform goodness-of-fit test for missing mechanisms 1 to 4.
Under the constraint $T\ind X$, Model \ref{miss::1} to Model \ref{miss::4} all have $10$ parameters and the observed data provides $11$ frequencies. Therefore,  the likelihood ratio test with an asymptotic $\chi^2(1)$ distribution can be used to test goodness-of-fit, i.e.,
$$
LR = 2 \sum\limits_{t, x, y = 0,1; m=0} N_{txy0} \log (N_{txy0} / \hat{N}_{txy0}) + 2 \sum\limits_{t,  y = 0,1; m=1} N_{t+y1} \log (N_{t+y1} / \hat{N}_{t+y1}) 
\stackrel{a}{\sim} \chi^2(1).
$$

Another interesting problem, raised by the Associate Editor, is testing the interaction of the treatment $T$ and the covariate $X$ on the outcome $Y$. 
The interaction is also called treatment heterogeneity, or effect modification.  
We can perform the likelihood ratio test for the interaction, which requires calculations of the likelihoods with and without the interaction of $T$ and $X$ on $Y.$
A more directly way is to compare the Bayesian posterior distributions of $CE_0$ and $CE_1$, or to find the credible interval of $CE_0 - CE_1$. If the credible interval of $CE_0 - CE_1$ does not contain $0$, we then find evidence of effect modification of $X$.

In our application in a randomized experiment, we will perform both the goodness-of-fit test and the effect modification test, under each missing data mechanism.


\section{Simulation Study}\label{sec::simulation}

In this section, we evaluate the finite sample performances of the likelihood-based and Bayesian inference for the missing mechanisms \ref{miss::1} to \ref{miss::4}, via a simulation study.
In order to mimic the real data analyzed in the next section, we assume that $T$ is completely randomized and thus $T\ind X$.
We generated $T\sim \text{Bernoulli}(0.5)$ and $X \sim \text{Bernoulli}(0.5)$.
Define $p_{y|tx} = P(Y=y\mid T=t, X=x)$, and we generated $Y$ according to the conditional distribution
$(p_{1|00}, p_{1|01}, p_{1|10}, p_{1|11}) = (0.2, 0.5, 0.8, 0.3)$ in all the cases.
We set the five missing data mechanisms to have the following parameters:

\begin{enumerate}[$M_1.~$]
\item $P(M=1\mid T=1, Y=1) = 0.3, P(M=1\mid T=1, Y=0) = 0.3, P(M=1\mid T=0, Y=1) = 0.4, P(M=1\mid T=0, Y=0)=0.7$;

\item $P(M=1\mid T=1, X=1) = 0.7, P(M=1\mid T=1, X=0) = 0.6, P(M=1\mid T=0, X=1) = 0.5, P(M=1\mid T=0, X=0)=0.3$;

\item $P(M=1\mid Y=1, X=1) = 0.3, P(M=1\mid Y=1, X=0) = 0.5, P(M=1\mid Y=0, X=1) = 0.3, P(M=1\mid Y=0, X=0)=0.8$;

\item $\text{logit}\{ P(M=1\mid T=t, X=x, Y=y) \}  = -1 + 1.4t - 0.5x + 0.8y$;

\item $\text{logit}\{ P(M=1\mid T=t, X=x, Y=y) \}  = -1 + 1.4t - x -0.5y + 0.5 tx + 0.3 ty - 0.6 xy - 0.2txy$.
\end{enumerate}

\begin{figure}
\begin{tabular}{p{\columnwidth}}
\includegraphics[width = \textwidth]{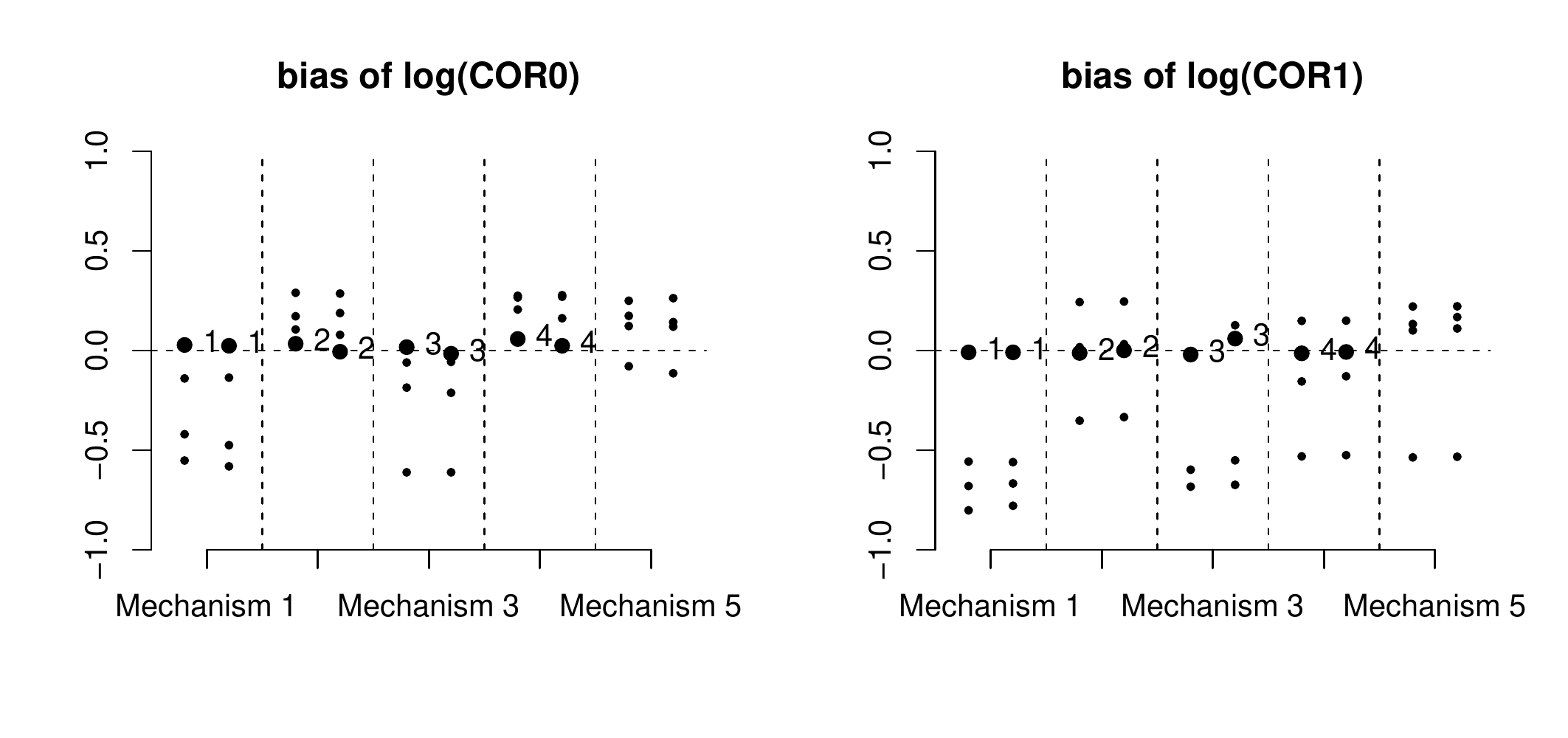}
\\
(a) Average biases. Each block corresponds to different missing data mechanism for generating the data. Within each block, there are two columns of points, corresponding to the biases of the MLEs and posterior medians. 
We label the average biases under the ``correct missing data mechanisms''.
\\
\includegraphics[width = \textwidth]{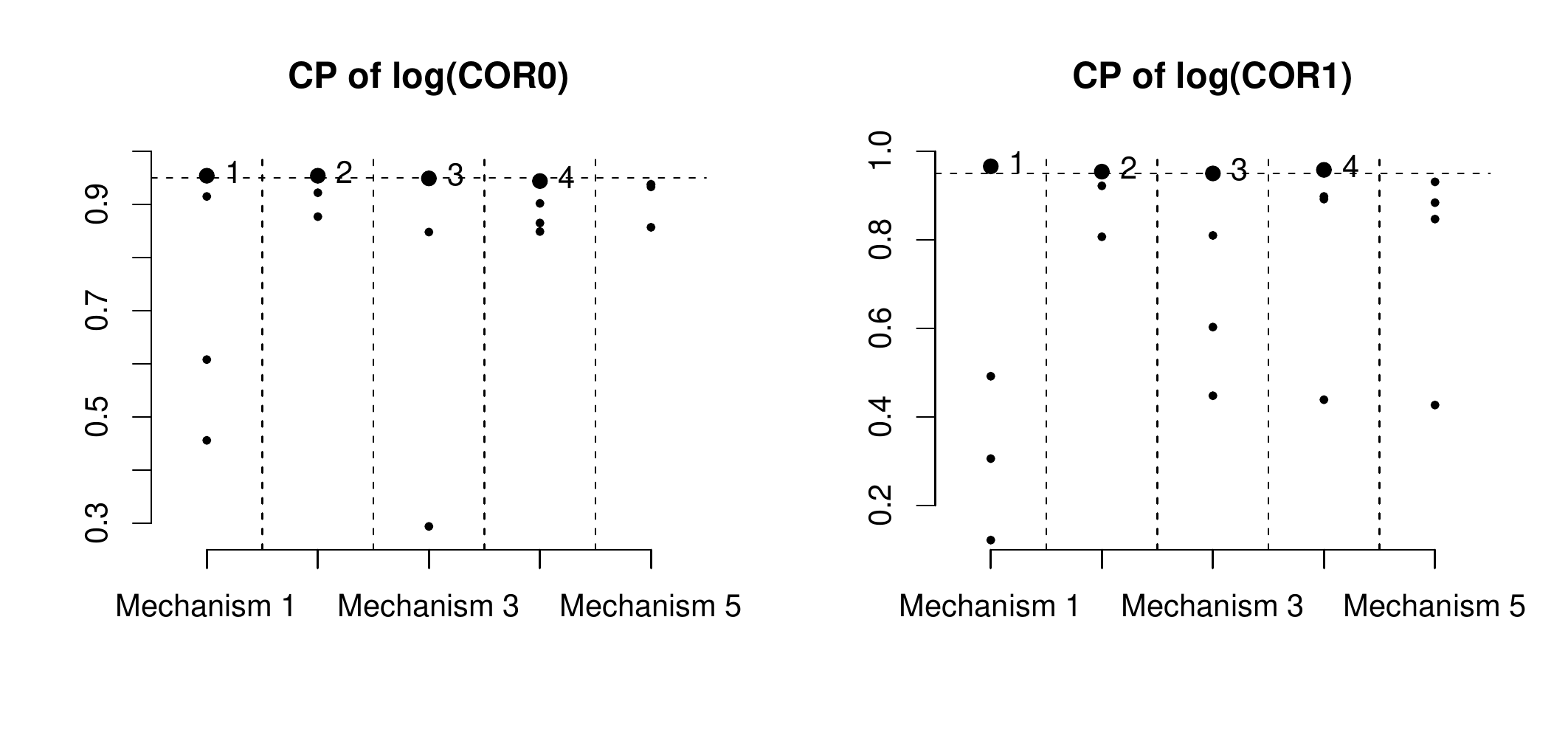}
\\
(b) Coverage proportions of the $95\%$ credible intervals.  
Each block corresponds to different missing data mechanisms, which shows the coverage proportions of the credible intervals under different models. 
We label the coverage proportions under the ``correct missing data mechanisms''.
\end{tabular}
\caption{Simulation Studies
} \label{fg::simulation}
\end{figure}

We apply the methods under missing mechanisms 1 to 4 to all of the five data sets. 
Thus we also show the sensitivities of our methods,
when the missing mechanisms are not correctly specified.
We use the EM algorithms to find the MLEs and use the Gibbs Samplers to find the posterior distributions of $ \log(COR_0)$ and $ \log(COR_1)$.
Using the Gibbs Samplers to obtain the Bayesian credible intervals is more direct than using the likelihood-based inference to obtain the confidence intervals.
The Gibbs samplers were run $10000$ times with burn-in after the $5000$-th iteration.
We did the simulation studies under sample sizes $500$ and $1000$, and the processes were repeated $1000$ times.
In Figure \ref{fg::simulation}, we show the simulation results under sample size $1000$, and detailed comparison of the results under sample sizes $500$ and $1000$ are presented in the web Appendix.

In Figure \ref{fg::simulation}(a), we show the average biases for $\log (COR_0)$ and $\log(COR_1)$ of both the MLEs and posterior medians.
For example, in the subfigure for ``bias of $\log(COR_0)$'', we divided the results into five blocks corresponding to the five data generating processes above. Within each block, there are two columns of points, corresponding to the biases of the MLEs on the left and posterior medians on the right. 
The results from MLEs and posterior medians are very similar to each other.
Clearly, if the missing data mechanisms are correctly specified, the average biases are very close to $0$. 
We label the average biases under the ``correct models'' for each data generating process, and all of them are very close to the horizontal zero line.
However, the average biases can be very arbitrary, if the missing data mechanisms are misspecified.
We have the same pattern in the subfigure for ``bias of $\log(COR_1)$''.

In Figure \ref{fg::simulation}(b), we show the coverage proportions of the $95\%$ credible intervals obtained from Bayesian posterior distributions of $\log(COR_0)$ and $\log(COR_1)$.
In the following, we will describe the subfigure for ``CP of $\log(COR_0)$'', and the same interpretation applies to the subfigure for ``CP of $\log(COR_1)$''.
Similar to the structures of the figures for the biases, we divide the results into five blocks, corresponding to five missing data mechanisms. We label the coverage proportions under the ``correct models'', which are very close to the nominal level $95\%$. But the coverage proportions can be extremely poor under model misspecifications.

\section{Application to a Randomized Clinical Trial}\label{sec::real}
\subsection{Background of the Motivating Example}
In this section, we will re-analyzed a randomized clinical trial [5]
using the newly proposed methods under different missing mechanisms.
We first briefly review the background of the illustrative clinical trial, and more details of the data can be found in the previous paper [5]. 
In this example, $T$ is the treatment assignment variable, with $T=1$ denoting the treatment (implantable cardiac defibrillator: ICD) and $T=0$ denoting the control. The outcome $Y$ is the death indicator, with $Y=1$ denoting ``dead'' and $Y=0$ denoting ``alive''.
Since there is no missing data in $Y$, it is direct to evaluate the causal effect of the treatment on the primary outcome of interest.
However, practitioners are also interested in evaluating the subgroup causal effects, stratified by the inducibility status. Let
$X$ denote the inducible indicator, with $X=1$ denoting ``inducible'' and $X=0$ denoting ``non-inducible''.
The covariate $X$ is obtained from the electro-physiological stimulation (EPS) testing. Since the EPS testing is invasive and not a prerequisite for enrollment in the study, $79.3\%$  of patients in the ICD arm have EPS records, while only $2.4\%$ of patients in the control arm have EPS records. Therefore, the missing data problem for the covariate $X$ is very severe.
The observed data can be summarized as the following counts: 
$N_{0000} = 4$, $ N_{0010} = 0$, $N_{1000} = 311$, $N_{1010} = 62$,
$N_{0100} = 6$, $ N_{0110} = 2$, $N_{1100} = 190$, $N_{1110} = 20$,
$N_{0+01} = 382$, $N_{0+11} = 95$, $N_{1+01} = 136$, and $N_{1+11} = 23$,
with the counts $N_{txy1}$ and $N_{t+y0}$ defined in Section 2.

Although the treatment assignment $T$ is independent of the inducibility status $X$ by randomization ($T\ind X$), the decision to conduct the test is associated with treatment assignment ($M\nind T$) and it may depend on the ``baseline'' factors that are associated with inducibility status ($M\nind X$) and the mortality ($M\nind Y$).

\subsection{Analysis of the Data, Model Criticism and Selection}

By randomization, the estimate of $\log(COR)$ for the population, $\log\left\{  \frac{ P(Y(1) = 1)P(Y(0) = 0)}{P(Y(1) = 0)P(Y(0) = 1)}  \right\}$, is $-0.235$ with the standard error $0.156$ based on the normal approximation, which is not significant at the level of $95\%$.
Since there is a zero cell count in the real data, the bounds of $\log(COR_x)$ obtained under missing mechanism \ref{miss::5} contain infinity and we do not present them here.

As discussed in Section \ref{sec::test}, randomization gives us one extra degree of freedom, and allows us to perform the likelihood ratio tests for goodness-of-fit test.
Table \ref{tb::mc} shows the results of the likelihood ratio tests, missing mechanisms \ref{miss::2} and \ref{miss::4} cannot be rejected, but missing mechanisms \ref{miss::1} and \ref{miss::3} are rejected. Complete randomization allows us to reject the ignorable missing data mechanism \ref{miss::1}, and one of the nonignorable missing data mechanism (mechanism 3).
The other two nonignorable missing mechanisms 2 and 4 fit the data very well, with mechanism 2 slightly better than mechanism 4.

Also, the following four plausible scientific assumptions are available based on the previous studies and some expert opinions [5]. 

\begin{assumption}\label{assume::2}
$P(X=0\mid T=t, M=1) \geq P(X=0\mid T=t, M=0)$ for $t=0,1.$
\end{assumption}
\begin{assumption}\label{assume::3}
$P(Y=1\mid T=0, X=1) \geq P(Y=1\mid T=0, X=0)$.
\end{assumption}
\begin{assumption}\label{assume::4}
$0.05\leq P(Y=1\mid T=0, X=x)\leq 0.50$ for $x=0,1$.
\end{assumption}
\begin{assumption}\label{assume::5}
$P(Y =1 \mid T =1, X=1) \leq P(Y =1\mid T =0, X=1).$
\end{assumption}

In Table \ref{tb::mc}, we also check whether these assumptions hold at the MLEs under each of the missing mechanism, where ``True'' denotes that an assumption is not rejected
and ``False'' denotes that it is rejected. 
The results are very interesting that the statistical findings are compatible with the clinical background and the expert opinions, because missing mechanisms \ref{miss::1} and \ref{miss::3} violate some of the expert opinions, and missing mechanisms \ref{miss::2} and \ref{miss::4} conform to all the expert opinions.
From Table \ref{tb::mc}, both missing mechanisms \ref{miss::2} and \ref{miss::4} fit the observed data very well and satisfy the scientific assumptions.
Based on the log likelihood and the scientific assumptions, missing mechanism \ref{miss::2} is chosen by us, which indicates that the missing data mechanism depends on the treatment assignment $T$ and the covariate $X$. 
Fortunately, the conclusion is not sensitive to the choice between missing mechanisms \ref{miss::2} and \ref{miss::4}, and we will only discuss the results under missing mechanism \ref{miss::2}.
Under missing mechanism \ref{miss::2}, the treatment is not significantly positively effective for neither the inducible subjects $(X=1)$ nor the non-inducible subjects ($X=0$) at the $95\%$ level.
We plot both the MLEs and posterior distributions of $\log(COR_1), \log(COR_0)$ and $\log(COR_1) - \log(COR_0)$ in Figure \ref{fg::model2-real}(a).
Although the posterior distributions of $\log(COR_1)$ and $ \log(COR_0)$ seems different in the left panel of Figure \ref{fg::model2-real}(a), the $95\%$ credible interval of $\log(COR_1) - \log(COR_0)$ contains $0$. Therefore, we conclude that the evidence of the effect modification of $X$ is not strong enough in this example.

\subsection{Comparison with the Original Analysis}

It would also be interesting to compare our analysis to the original analysis [5], which used causal risk ratios as the causal measures of interest.
Under missing mechanism 2, the MLEs are $\widehat{CRR}_1 = 0.301$ and $ \widehat{CRR}_0 = 1.279$. The posterior median of $CRR_1$ is $0.303$ with a $95\%$ credible interval $[0.140, 1.176]$, and the posterior median of $CRR_0$ is $1.551$ with a $95\%$ credible interval $[0.593, 227.278]$. The wide interval for $CRR_0$ is due to the large proportion of missing data and heavy-tailedness of the posterior distribution. 
However, the original analysis did not provide us with point estimators, and their bounds for $CRR_1$ are $[0.05, 46.01]$, and their bounds for $CRR_0$ are $[0.21, 1425.58]$. 
Both our MLEs and posterior medians are within the bounds obtained in the previous study.
The bounds obtained by the previous analysis is much wider than ours, and they can only be sharpened with the help of the expert opinions.

\subsection{Sensitivity Analysis}

\begin{figure}[htb]
\begin{tabular}{p{\columnwidth}}
\includegraphics[width = \textwidth]{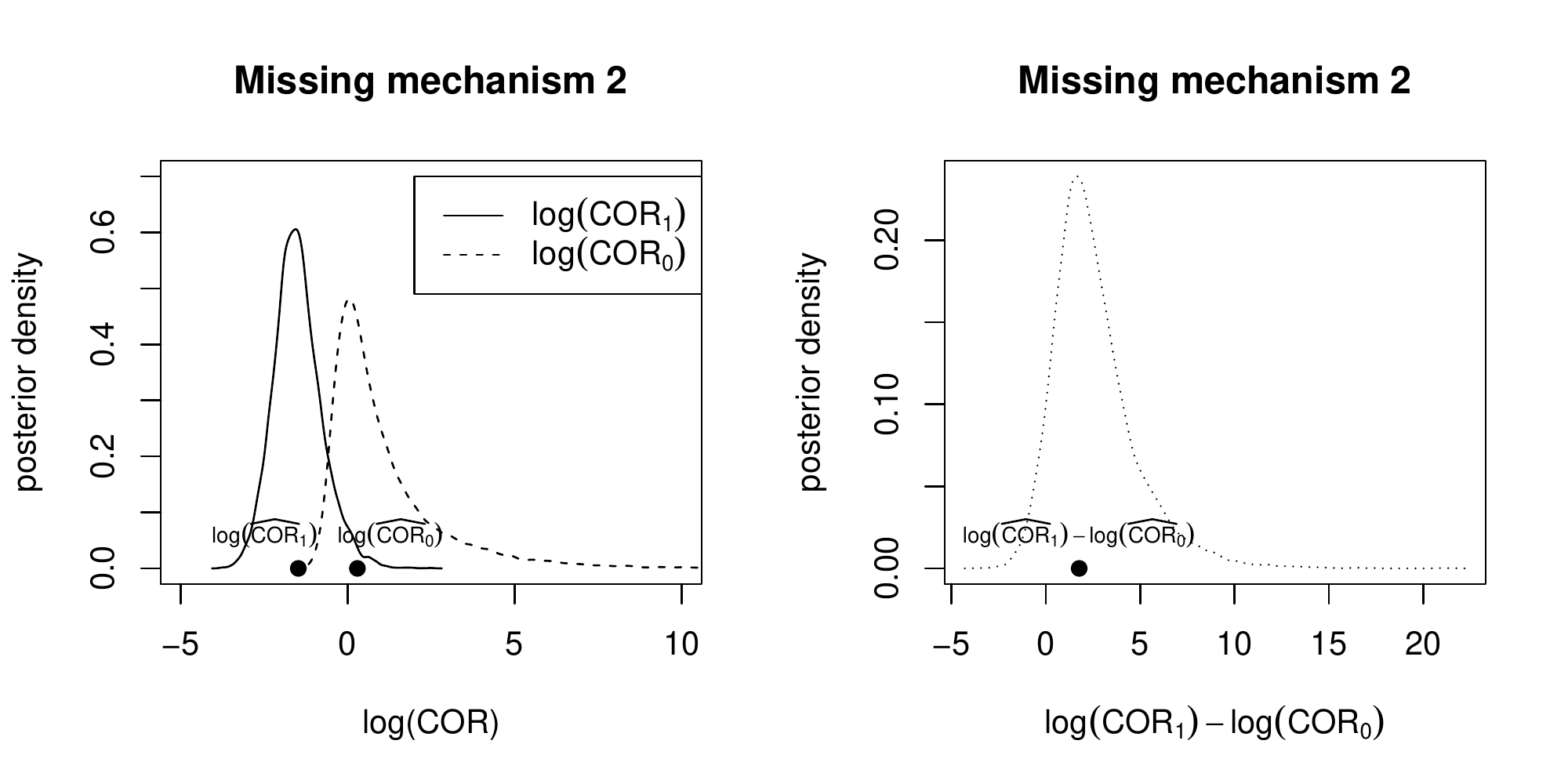}\\
(a) MLEs and posterior distributions under missing mechanism 2.
\\
\includegraphics[width=\textwidth]{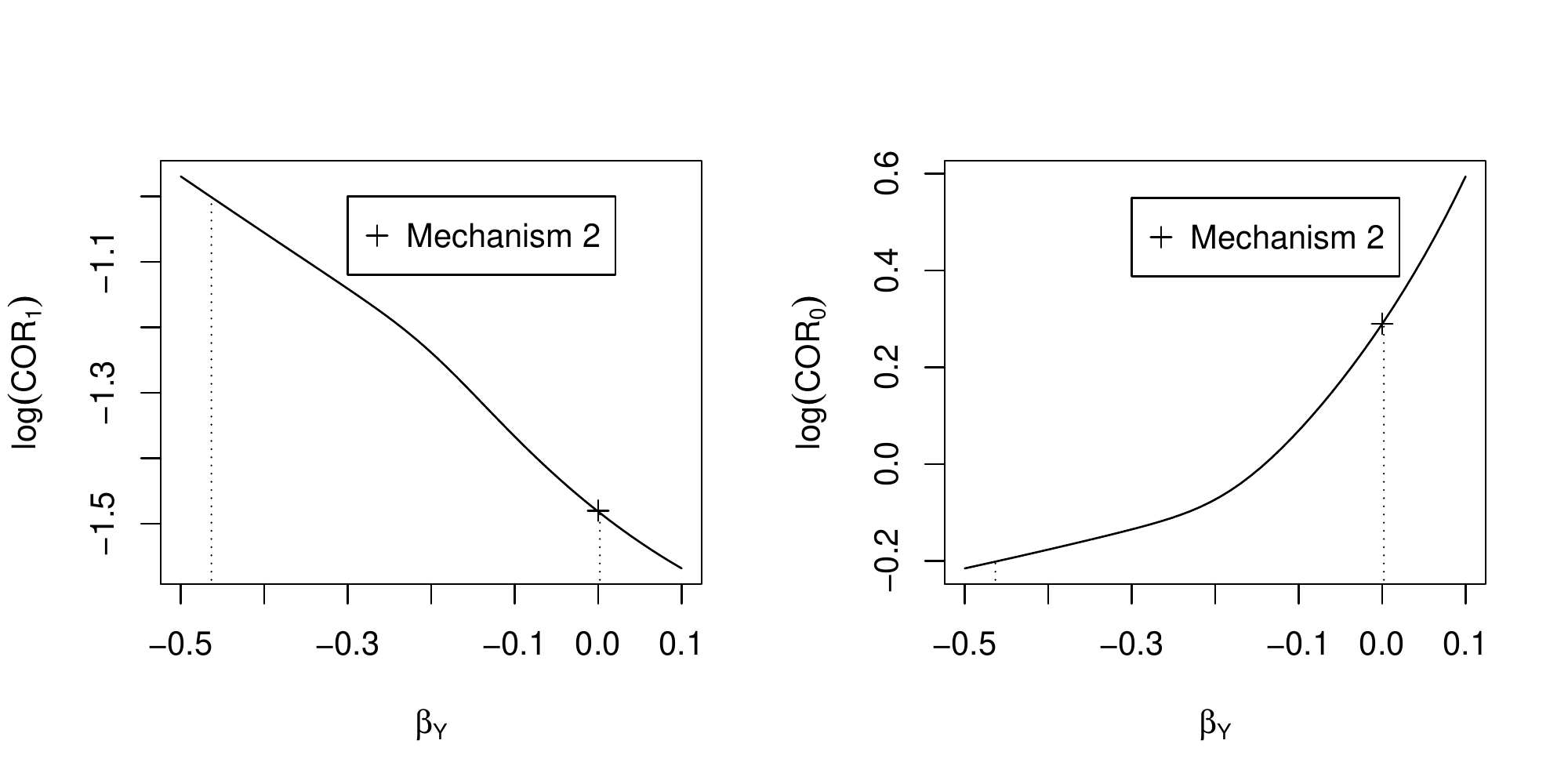}\\
(b) Sensitivity analysis for missing mechanism 2. The parameter $\beta_Y$ in (\ref{eq::sensi}) is the sensitivity parameter, and
the ``feasible'' regions are within the dotted lines. Here, ``feasible'' means that the MLEs within these regions are compatible with the expert opinions.
\end{tabular}
\caption{Analysis of the data under missing mechanism 2}\label{fg::model2-real}
\end{figure}

The saturated model for the missing data is
$$
\text{logit} \{ P(M=0\mid T=t, X=x, Y=y) \} = \beta_0 + \beta_T t + \beta_X x + \beta_Y y + \beta_{TX} tx + \beta_{TY} ty + \beta_{XY} xy + \beta_{TXY} txy.
$$
Each of the missing mechanisms from \ref{miss::1} to \ref{miss::4} restricts $4$ coefficients $\beta$'s to be $0$. For example, missing mechanism 1 restricts $\beta_X = \beta_{TX} = \beta_{XY} = \beta_{TXY} = 0$, and analogous results hold for other missing mechanisms.
We treat missing mechanism \ref{miss::2} as the benchmark for our sensitivity analysis and discuss the following model which allows $M$ to depend on $T,X$ and $Y$:
\begin{eqnarray}
\text{logit} \{ P(M=0\mid T=t, X=x, Y=y) \} = \beta_0 + \beta_T t + \beta_X x + \beta_{TX} tx  +  \beta_Y y,   \label{eq::sensi}
\end{eqnarray}
where $\beta_Y$ is the sensitivity parameter, and $\beta_Y=0$ corresponds to missing mechanism \ref{miss::2}.
Fixing $\beta_Y$ at different values and obtaining MLEs for Model (\ref{eq::sensi}), the sensitivities of $\log(COR_1)$ and $\log(COR_0)$ are shown in Figure \ref{fg::model2-real}(b), where the ``feasible'' regions are within the dotted lines. Here, ``feasible'' means that 
the MLEs of 
Model (\ref{eq::sensi})  are compatible with the expert opinions, when $\beta_Y$ is within these regions.
Positive $\beta_Y$'s in Figure \ref{fg::model2-real}(b) are implausible under Model (\ref{eq::sensi}), and $\log(COR_1)$ are less sensitive than $\log(COR_0)$, since $\log(COR_1)$ does not change sign within the ``feasible" region.

\begin{table}[ht]
\begin{center}
\caption{Model Comparison. The ``log likelihood'' column shows the log likelihoods evaluated at the MLEs.
The ``$p$-value of LRT'' column shows the $p$-values of the likelihood ratio tests for goodness-of-fit. The ``A.\ref{assume::2}'' to ``A.\ref{assume::5}'' columns check whether Assumptions \ref{assume::2} to \ref{assume::5} hold for the MLEs under mechanisms \ref{miss::1} to \ref{miss::4}, with ``True'' if one assumption is not rejected and ``False'' if one assumption is rejected.
}
\label{tb::mc}
\begin{tabular}{rrrrrrr}
  \hline
Mechanism& log likelihood & $p$-value of LRT & A.\ref{assume::2}& A.\ref{assume::3}& A.\ref{assume::4}& A.\ref{assume::5} \\
  \hline
  Mechanism \ref{miss::1}& -2202.654  & 0.017       &False&True&True&True\\
Mechanism \ref{miss::2} &  -2200.452 & 0.248      &True&True&True&True \\
Mechanism \ref{miss::3}& -2503.779  & $<$0.001 &False&False&False&False\\
Mechanism \ref{miss::4}& -2200.584 & 0.206      &True&True&True&True\\
   \hline
\end{tabular}
\end{center}
\end{table}

\section{Application to a Job-Training Program Data}
\label{sec::job-train}

In this section, we use a data set from a job-training program [18]
to illustrate the potential applications of our methods to observational studies. The data is available in the R package ``Matching'' [19], 
and more detailed descriptions of the data can be found in the previous papers [18, 19]. 
The data set contains $445$ observations of males, with the treatment $T$ as the indicator for receiving the job-training program, and the outcome $Y$ as the indicator for unemployment in 1978. The pretreatment covariates in the original data set are age, education, ethnicity, marriage status, historical employment status, and the indicator for a high school degree.
Although the original data is a completely randomized experiment, we found some evidence of imbalance in the pretreatment covariates. 
All covariates are balanced between the treatment and control groups except the binary ``no degree'' variable, and the balance checking result is shown in the web Appendix.
In order to adjust for the imbalance of the covariate $X=$ ``no degree'', we assume that the experiment is randomized conditionally on the covariate $X$. 
Although this data set is not from an observational study, it has the same nature of an observational study. And our analysis is under Assumption 2 instead of Assumption 3.

Another feature of the data is that all the values of the covariate $X$ are observed, therefore we can identify the subgroup causal effects and we know the ``true'' MLEs from the complete data. Recall the definition of $p_{y|tx} = P(Y=y\mid T=t, X=x)$ in Section 6, and we find their MLEs as $  (\widehat{p}_{1|11}, \widehat{p}_{1|01},  \widehat{p}_{1|10}, \widehat{p}_{1|00} ) = (0.260, 0.355, 0.204, 0.349).$
Thus, the MLEs for the subgroup causal effects in term of the risk difference are $\widehat{CRD}_0 = -0.095$ ($p$-value = 0.058) and $\widehat{CRD}_1 = -0.145$ ($p$-value = 0.111). And the test for effect modification has $p$-value 0.632, providing very weak evidence of effect modification.

To illustrate our methods in observational studies, we artificially create missing values in the covariate $X$ according to different underlying mechanisms, and try to recover the ``true'' MLEs using our methods.
The following four cases correspond to four missing data mechanisms.

\begin{enumerate}
[$M_1.~$]
\item We create missing data according to $P(M=1\mid T=1, Y=1)=0.4, P(M=1\mid T=0, Y=1 )= 0.3, P(M=1\mid T=1, Y=0) = P(M=1\mid T=0, Y= 0) = 0.2$;

\item We create missing data according to $P(M=1\mid T=1, X= 1)=0.5, P(M=1\mid T=0, X=1) = 0.3, P(M=1\mid T=1, X=0) = P(M=1\mid T=0, X=0) = 0.2$;

\item We create missing data according to $P(M=1\mid X=1, Y=1) = 0.6, P(M=1\mid X=1, Y=0) = P(M=1\mid X=0, Y=1) = P(M=1\mid X=0, Y=0) = 0.1$;

\item We create missing data according to $\text{logit}\{ P(M=1\mid T=t, X=x, Y=y) \}  = - 1 + t -x+y$.

\end{enumerate}

We then fit each of the generated data set by MLEs under these missing mechanisms 1 to 4, corresponding to estimation methods $h=1,2,3,4.$
Denote $ (\widehat{p}_{1|11}^h, \widehat{p}_{1|01}^h,  \widehat{p}_{1|10}^h, \widehat{p}_{1|00}^h ) $ as the MLEs from method based on the missing mechanism $h$. We use $ RMSE(h) =\sqrt{   \sum_{t,x}  (  \widehat{p}_{1|tx} - \widehat{p}_{1|tx}^h )^2 }$ as the criterion to evaluate the results from different estimation methods.
The rows of Table \ref{tb::job-train} correspond to different missing data mechanisms for generating the data, and the columns correspond to different estimation methods.
Table \ref{tb::job-train} shows that when the missing data mechanisms are correctly specified, we can recover the ``true'' MLEs very well, that is, the diagonal elements in Table \ref{tb::job-train} is the smallest in each row.
 However, with misspecified missing data mechanisms, the behaviors of the MLEs may be very arbitrary and far from the ``true'' MLEs.
This illustrative example demonstrates the importance of the specification of the missing data mechanism in the analysis of observational studies, when some key covariates are missing.

\begin{table}[ht]
\caption{Analysis of the Job-Training Program. The rows correspond to different missing data generating mechanisms, and the columns correspond to different estimation methods.}
\label{tb::job-train}
\centering
\begin{tabular}{ccccc}
\hline
&\multicolumn{4}{c}{Mechanism used in estimation}\\
\cline{2-5}
Missing data mechanism & 1 & 2 &  3 & 4\\
 \hline
Mechanism 1& \textbf{0.014} & 0.118 & 0.270 & 0.110\\
Mechanism 2& 0.009 & \textbf{0.007} & 0.075 & 0.010\\
Mechanism 3&0.355 & 0.185 &  \textbf{0.031}& 0.326\\
Mechanism 4&0.150 & 0.055& 0.296 &\textbf{0.031}\\
\hline
\end{tabular}
\end{table}


\section{Discussion}\label{sec::discussion}

Randomized experiments are regarded as the gold standard for causal evaluations.
In the cases with a nonignorable missing covariate,
however, the evaluation of the subgroup causal effects
conditional on the covariate is very challenging. 
Without making untestable assumptions, we can only get bounds of the subgroup causal effects, which may be too wide to be useful. 
We show that they are pointly identifiable, and we perform both likelihood-based and Bayesian inference, under some model assumptions and restrictions.
An interesting application of our proposed models and methods to a randomized experiment shows that one of the nonignorable missing data mechanisms is more appropriate than the ignorable missing data mechanism, with a higher likelihood in the former model than the latter one.
Expert opinions help us to verify our statistical findings, because the results from the chosen nonignorable missing data mechanism are consistent with the expert opinions, while the results from the ignorable missing data mechanism are not.
Since many fundamental assumptions in causal inference are not directly testable, utilizing expert opinions to guide our practice is very valuable and should be tirelessly emphasized.

There are several issues beyond the scope of this paper, and our discussion below benefits a lot from the Associate Editor and a reviewer's comments. 
First, we discussed the identifiability of subgroup causal effects when the outcome $Y$ is categorical.
However, continuous outcomes are also very common in practice.
One approximate approach is to categorize the continuous outcome $Y$. Another approach for the missing mechanisms 1 and 2 is to dichotomize the continuous $Y$ as $I(Y>y)$ for each observed value $y$,
and we can identify the average causal effect of $Z$ on $I_y$, which is also the distributional causal effect of $Z$ on $Y$ at point $y$. 
The subgroup causal effects are identifiable, once the subgroup distributional causal effect is identified [20].
When the missing data mechanism depends on the outcome (mechanisms 3 and 4), 
we need to simultaneously model the outcomes and the nonignorable missing data mechanism.

Second, it is possible that the covariates have high dimensions. 
Identifying subgroup causal effects defined by high dimensional categorical covariates is difficult even if there is no missing data, since the observations within each subgroup can be very sparse in finite samples. 
If the subgroups are defined by only one categorical variable subject to missing but there are a large number of other categorical and/or continuous variables which are completely observed or missing at random,
our identifiability results also hold, but some approaches for high dimensional covariates should be used for the estimation. In the cases where there are more than one covariate subject to nonignorable missing, more complicated missing mechanisems must be introduced. 
Other approaches such as multiple imputation [21]
and jointly modeling the distribution of $X$ and the missing data mechanism [17] 
may be used to treat the high dimensional missing covariates.

Third, we discussed the case in which only the covariate may be missing not at random. In many applications, both the covariate and the outcome may be missing not at random. Then we must describe the missing mechanisms for both of them, and the identifiability and estimation would be more complicated.

All the topics mentioned above are of great interest both theoretically and practically. Although they are beyond our current study, we will go on our research in this area.

\section*{Acknowledgement}
We would like to thank both the Associate Editor and the reviewer for their very insightful and constructive comments, which lead to a significant improvement of our paper.
This research was supported by NSFC (11171365, 11021463, 10931002).


\section*{References}

\begin{itemize}

\item[1.]
Egleston BL, Scharfstein DO, and MacKenzie E. On estimation of the survivor average causal effect in observational studies when important confounders are missing due to death. {\it Biometrics} 2009; {\bfseries 65}(2):497-504.

\item[2.] 
Frangakis CE, Rubin DB, An MW, and MacKenzie E. Principal stratification designs to estimate input data missing due to death (with discussion). {\it Biometrics} 2007; {\bfseries 63}(3):641-649.

\item[3.]
Horowitz JL and Manski CF. Nonparametric analysis of randomized experiments with missing covariate and outcome data (with discussion). {\it Journal of the American Statistical Association} 2000; {\bfseries 95}(449):77-84.

\item[4.] 
Rathouz PJ. Identifiability assumptions for missing covariate data in failure time regression models. {\it Biostatistics} 2007; {\bfseries 8}(2):345-356.

\item[5.]
Scharfstein D, Onicescu G, Goodman S, and Whitaker R. Analysis of subgroup effects in randomized trials when subgroup membership is missing: application to the second multicenter automatic defibrillator intervention trial. {\it Journal of the Royal Statistical Society: Series C (Applied Statistics)} 2011; {\bfseries 60}(4):607-617.

\item[6.] 
Yan W, Hu Y, and Geng Z. Identifiability of causal effects for binary variables with baseline data missing due to death. {\it Biometrics} 2012; {\bfseries 68}(1):121-128.

\item[7.]
Chen H, Geng Z, and Zhou XH. Identifiability and estimation of causal effects in randomized trials with noncompliance and completely nonignorable missing data (with discussion). {\it Biometrics} 2009; {\bfseries 65}(3):675-682.

\item[8.] 
Imai K. Statistical analysis of randomized experiments with non-ignorable missing binary outcomes: an application to a voting experiment. {\it Journal of the Royal Statistical Society: Series C (Applied Statistics)} 2009; {\bfseries 58}(1):83-104.

\item[9.]
Ma WQ, Geng Z, and Hu YH. Identification of graphical models for nonignorable nonresponse of binary outcomes in longitudinal studies. {\it Journal of Multivariate Analysis} 2003; {\bfseries 87}(1):24-45.

\item[10.]
Egleston BL and Wong YN. Sensitivity analysis to investigate the impact of a missing covariate on survival analyses using cancer registry data. {\it Statistics in Medicine} 2009; {\bfseries 28}(10):1498-1511.

\item[11.]
Frangakis CE and Rubin DB. Addressing complications of intention-to-treat analysis in the combined presence of all-or-none treatment-noncompliance and subsequent missing outcomes. {\it Biometrika} 1999; {\bfseries 86}(2):365-379.

\item[12.]
Little RJ and Zhang N. Subsample ignorable likelihood for regression analysis with missing data. {\it Journal of the Royal Statistical Society: Series C (Applied Statistics)} 2011; {\bfseries 60}(4):591-605.

\item[13.]
Rubin DB. Comment on ``Randomization analysis of experimental data: The Fisher randomization test''. {\it Journal of the American Statistical Association} 1980; {\bfseries 75}(371):591-593.

\item[14.] 
Rubin DB. Comment on ``Statistics and causal inference'': Which ifs have causal answers. {\it Journal of the American Statistical Association} 1986; {\bfseries 81}(396):961-962.

\item[15.]
Rubin DB. Bayesian inference for causal effects: The role of randomization. {\it The Annals of Statistics} 1978; {\bfseries 6}(1):34-58.

\item[16.]
Jin H and Rubin DB. Principal stratification for causal inference with extended partial compliance. {\it Journal of the American Statistical Association} 2008; {\bfseries 103}(481):101-111.

\item[17.]
Little RJ and Rubin DB. {\it Statistical Analysis with Missing Data}. Wiley New York, 2002.

\item[18.]
Dehejia RH and Wahba S. Causal effects in nonexperimental studies: Reevaluating the evaluation of training programs. {\it Journal of the American statistical Association} 1999; {\bfseries 94}(448):1053-1062.

\item[19.] 
Sekhon JS. Multivariate and propensity score matching software with automated balance optimization: The Matching package for R. {\it Journal of Statistical Software} 2011; {\bfseries 42}(7):1-52.

\item[20.]
Ding P, Geng Z, Yan W, and Zhou XH. Identifiability and estimation of causal effects by principal stratification with outcomes truncated by death. {\it Journal of the American Statistical Association} 2011; {\bfseries 106}(496):1578-1591.

\item[21.]
Rubin DB. {\it Multiple Imputation for Nonresponse in Surveys}. John Wiley \& Sons: New York, 1987.

\end{itemize}

\clearpage 

\newpage 

\begin{flushleft}
{\Large \bf Supporting Materials for ``Identifiability of Subgroup Causal Effects in Randomized Experiments with Nonignorable Missing Covariates'' by Peng Ding and Zhi Geng}\\

\end{flushleft}

\section{Computational Details}

\subsection{Missing Mechanism 1}
Denote $P^{(j)}(T=t, X=x, Y=y, M=m) = P^{(j)}(X=x) P^{(j)}(T=t\mid X=x)  P^{(j)}(Y=y\mid T=t, X=x) P^{(j)}(M= m \mid T=t, X=x)$ as the joint distribution of $(T,X,Y,M)$ in the $j$-th iteration. Define
\begin{eqnarray*}
 p_{x|ty1}^{(j)} &=& P^{(j)}(X=x\mid T=t, Y=y, M=1) \\
 &=& \frac{ P^{(j)}(T=t, X=x, Y=y, M=1) }  
 {  \sum_{x'=0,1} P^{(j)}(T=t, X=x', Y=y, M=1)  }.
\end{eqnarray*}
The EM algorithm iterates between the following two steps:
\begin{itemize}
\item
E-step: The sufficient statistics are imputed as $N^{(j)}_{txy0} = N_{txy0}$ and $N^{(j)}_{txy1} = N_{t+y1} p_{x|ty1}^{(j)} $;

\item
M-step: 
The joint distribution is updated by 
$P^{(j+1)}(T=t, X=x, Y=y, M=m) = \frac{N^{(j)}_{+x++}}{N^{(j)}_{++++}}   \frac{N^{(j)}_{tx++}}{N^{(j)}_{+x++}} 
 \frac{N^{(j)}_{txy+}}{N^{(j)}_{tx++}}  \frac{N^{(j)}_{tx+m}}{N^{(j)}_{tx++}}$.

\end{itemize}

%
\noindent The Gibbs sampler iterates between the following two steps:
\begin{itemize}
\item
Imputation step: 
We let $N^{(j)}_{txy0} = N_{txy0}$ and draw $N^{(j)}_{txy1} \sim \text{Binomial}(N_{t+y1}, p_{x|ty1}^{(j)}) $;

\item
Posterior step:
Draw 
$ P^{(j+1)}(X=1) \sim \text{Beta} (\alpha_X + N^{(j)}_{+1++}, \beta_X +N^{(j)}_{+0++} )$,
$P^{(j+1)}(T=1\mid X=x)\sim \text{Beta}( \alpha_T^x + N^{(j)}_{1x++}, \beta_{T}^x + N^{(j)}_{0x++}  )$,
$ P^{(j+1)}(Y=1\mid T=t, X=x) \sim \text{Beta} (\alpha_{Y}^{tx} + N^{(j)}_{tx1+} , \beta_{Y}^{tx} + N^{(j)}_{tx0+} )$,
$  P^{(j+1)}(M=1\mid T=t, X=x) \sim \text{Beta} (\alpha_{M}^{tx} + N^{(j)}_{tx+1} , \beta_{M}^{tx} + N^{(j)}_{tx+0} )$,
where $\alpha_X, \alpha_T^x, \beta_T^x,  \beta_X, \alpha_Y^{tx}, \beta_{Y}^{tx}, \alpha_M^{tx}$, and $ \beta_M^{tx}$ are parameters for the Beta priors of the probability parameters.
\end{itemize}

%

\subsection{Missing Mechanism 2}
Denote $P^{(j)}(T=t, X=x, Y=y, M=m) = P^{(j)}(X=x) P^{(j)}(T=t\mid X=x)  P^{(j)}(Y=y\mid T=t, X=x) P^{(j)}(M= m \mid T=t, Y=y)$ as the joint distribution of $(T,X,Y,M)$ in the $j$-th iteration. Define
\begin{eqnarray*}
 p_{x|ty1}^{(j)} &=& P^{(j)}(X=x\mid T=t, Y=y, M=1) \\
 &=& \frac{ P^{(j)}(T=t, X=x, Y=y, M=1) }  
 {  \sum_{x'=0,1} P^{(j)}(T=t, X=x', Y=y, M=1)  }.
\end{eqnarray*}

The EM algorithm iterates between the following two steps:
\begin{itemize}
\item
E-step: The sufficient statistics are imputed as $N^{(j)}_{txy0} = N_{txy0}$ and $N^{(j)}_{txy1} = N_{t+y1} p_{x|ty1}^{(j)} $;

\item
M-step: 
The joint distribution is updated by 
$P^{(j+1)}(T=t, X=x, Y=y, M=m) =  \frac{N^{(j)}_{+x++}}{N^{(j)}_{++++}}  \frac{N^{(j)}_{tx++}}{N^{(j)}_{+x++}}  \frac{N^{(j)}_{txy+} }{N^{(j)}_{tx++} }  \frac{N^{(j)}_{t+ym}}{N^{(j)}_{t+y+}}$.

\end{itemize}

\noindent The Gibbs sampler iterates between the following two steps:
\begin{itemize}
\item
Imputation step: 
We let $N^{(j)}_{txy0} = N_{txy0}$ and draw $N^{(j)}_{txy1} \sim \text{Binomial}(N_{t+y1}, p_{x|ty1}^{(j)}) $;

\item
Posterior step:
Draw 
$
 P^{(j+1)}(X=1) \sim \text{Beta} (\alpha_X + N^{(j)}_{+1++}, \beta_X +N^{(j)}_{+0++} ),
 P^{(j+1)}(T=1\mid X=x)\sim \text{Beta}( \alpha_T^x + N^{(j)}_{1x++}, \beta_{T}^x + N^{(j)}_{0x++}  ),
 P^{(j+1)}(Y=1\mid T=t, X=x) \sim \text{Beta} (\alpha_{Y}^{tx} + N^{(j)}_{tx1+} , \beta_{Y}^{tx} + N^{(j)}_{tx0+} ),
  P^{(j+1)}(M=1\mid T=t, Y=y) \sim \text{Beta} (\alpha_{M}^{ty} + N^{(j)}_{t+y1} , \beta_{M}^{ty} + N^{(j)}_{t+y0} ),$ where $ \alpha_X, \beta_X, \alpha_T^x, \beta_T^x, \alpha_Y^{tx}, \beta_{Y}^{tx}, \alpha_M^{ty}$, and $ \beta_M^{ty}$ are parameters for the Beta priors.
\end{itemize}

%
%

\subsection{Missing Mechanism 3}  

Denote $P^{(j)}(T=t, X=x, Y=y, M=m) = P^{(j)}(X=x) P^{(j)}(T=t\mid X=x)  P^{(j)}(Y=y\mid T=t, X=x) P^{(j)}(M= m \mid X=x, Y=y)$ as the joint distribution of $(T,X,Y,M)$ in the $j$-th iteration. Define
\begin{eqnarray*}
 p_{x|ty1}^{(j)} &=& P^{(j)}(X=x\mid T=t, Y=y, M=1) \\
 &=& \frac{ P^{(j)}(T=t, X=x, Y=y, M=1) }  
 {  \sum_{x'=0,1} P^{(j)}(T=t, X=x', Y=y, M=1)  }.
\end{eqnarray*}

The EM algorithm iterates between the following two steps:
\begin{itemize}
\item
E-step: The sufficient statistics are imputed as $N^{(j)}_{txy0} = N_{txy0}$ and $N^{(j)}_{txy1} = N_{t+y1} p_{x|ty1}^{(j)} $;

\item
M-step: 
The joint distribution is updated by 
$P^{(j+1)}(T=t, X=x, Y=y, M=m) =  \frac{N^{(j)}_{+x++}}{N^{(j)}_{++++}} \frac{N^{(j)}_{tx++}}{N^{(j)}_{+x++}}
 \frac{N^{(j)}_{txy+}}{N^{(j)}_{tx++}}  \frac{N^{(j)}_{+xym}}{N^{(j)}_{+xy+}}$.

\end{itemize}

\noindent The Gibbs sampler iterates between the following two steps:
\begin{itemize}
\item
Imputation step: 
We let $N^{(j)}_{txy0} = N_{txy0}$ and draw $N^{(j)}_{txy1} \sim \text{Binomial}(N_{t+y1}, p_{x|ty1}^{(j)}) $;

\item
Posterior step:
Draw 
$
 P^{(j+1)}(X=1) \sim \text{Beta} (\alpha_X + N^{(j)}_{+1++}, \beta_X +N^{(j)}_{+0++} ),
 P^{(j+1)}(T=1\mid X=x)\sim \text{Beta}( \alpha_T^x + N^{(j)}_{1x++}, \beta_{T}^x + N^{(j)}_{0x++}  ),
 P^{(j+1)}(Y=1\mid T=t, X=x) \sim \text{Beta} (\alpha_{Y}^{tx} + N^{(j)}_{tx1+} , \beta_{Y}^{tx} + N^{(j)}_{tx0+} ),
  P^{(j+1)}(M=1\mid  X=x, Y=y) \sim \text{Beta} (\alpha_{M}^{xy} + N^{(j)}_{+xy1} , \beta_{M}^{xy} + N^{(j)}_{+xy0} ),$ where $\alpha_T, \beta_T, \alpha_X, \beta_X, \alpha_Y^{tx}, \beta_{Y}^{tx}, \alpha_M^{xy}$, and $ \beta_M^{xy}$ are parameters for the Beta priors.
\end{itemize}

%

\subsection{Missing Mechanism 4}

Denote $P^{(j)}(T=t, X=x, Y=y, M=m) = P^{(j)}(X=x) P^{(j)}(T=t\mid X=x)  P^{(j)}(Y=y\mid T=t, X=x) P^{(j)}(M= m \mid T= t, X=x, Y=y)$ as the joint distribution of $(T,X,Y,M)$ in the $j$-th iteration. Define
\begin{eqnarray*}
&& p_{x|ty1}^{(j)} \\
 &=& P^{(j)}(X=x\mid T=t, Y=y, M=1) \\
 &=& \frac{ P^{(j)}(T=t, X=x, Y=y, M=1) }  
 {  \sum_{x'=0,1} P^{(j)}(T=t, X=x', Y=y, M=1)  }\\
 &=& \frac{   P^{(j)}(X=x) P^{(j)}(T=t\mid X=x)  P^{(j)}(Y=y\mid T=t, X=x) \text{expit}(\beta^{(j)}_0 + \beta^{(j)}_T t + \beta^{(j)}_X x + \beta^{(j)}_Y y)    }
 { \sum_{x'=0,1}  P^{(j)}(X=x) P^{(j)}(T=t\mid X=x)  P^{(j)}(Y=y\mid T=t, X=x) \text{expit}(\beta^{(j)}_0 + \beta^{(j)}_T t + \beta^{(j)}_X x' + \beta^{(j)}_Y y)   },
\end{eqnarray*}
where $\text{expit}(a) = 1/(1 + e^{-a})$ and $(\beta^{(j)}_0,  \beta^{(j)}_T,  \beta^{(j)}_X, \beta^{(j)}_Y)$ are the parameters of the missing data mechanism at the $j$-th iteration.

The EM algorithm iterates between the following two steps:
\begin{itemize}
\item
E-step: The sufficient statistics are imputed as $N^{(j)}_{txy0} = N_{txy0}$ and $N^{(j)}_{txy1} = N_{t+y1} p_{x|ty1}^{(j)} $;

\item
M-step: 
The joint distribution is updated by 
$P^{(j+1)}(T=t, X=x, Y=y) =  \frac{N^{(j)}_{+x++}}{N^{(j)}_{++++}}   \frac{N^{(j)}_{tx++}}{N^{(j)}_{+x++}}    \frac{N^{(j)}_{txy+}}{N^{(j)}_{tx++}} $ and $(\beta^{(j+1)}_0,  \beta^{(j+1)}_T,  \beta^{(j+1)}_X, \beta^{(j+1)}_Y)$ are obtained by a Logistic regression of $\{ (N^{(j)}_{txy1}, N^{(j)}_{txy0}) \}$ on $(T=t, X=x, Y=y)$ where $t,x,y=0,1.$

\end{itemize}

\noindent The Gibbs sampler iterates between the following two steps:
\begin{itemize}
\item
Imputation step: 
We let $N^{(j)}_{txy0} = N_{txy0}$ and draw $N^{(j)}_{txy1} \sim \text{Binomial}(N_{t+y1}, p_{x|ty1}^{(j)}) $;

\item
Posterior step:
Draw 
$
 P^{(j+1)}(X=1) \sim \text{Beta} (\alpha_X + N^{(j)}_{+1++}, \beta_X +N^{(j)}_{+0++} ),
 P^{(j+1)}(T=1\mid X=x)\sim \text{Beta}( \alpha_T^x + N^{(j)}_{1x++}, \beta_{T}^x + N^{(j)}_{0x++}  ),
 P^{(j+1)}(Y=1\mid T=t, X=x) \sim \text{Beta} (\alpha_{Y}^{tx} + N^{(j)}_{tx1+} , \beta_{Y}^{tx} + N^{(j)}_{tx0+} )$,
and $(\beta^{(j+1)}_0,  \beta^{(j+1)}_T,  \beta^{(j+1)}_X, \beta^{(j+1)}_Y)$ are drawn by the Metropolis-Hastings algorithm, where $\alpha_T, \beta_T, \alpha_X, \beta_X, \alpha_Y^{tx}, \beta_{Y}^{tx}$ are parameters for the Beta priors of the probability parameters and the priors for $(\beta_0,  \beta_T,  \beta_X, \beta_Y)$ are flat.
\end{itemize}

\section{Proofs of the Theorems}
\setcounter{equation}{0}
\renewcommand {\theequation} {A.\arabic{equation}}

In the proofs of theorems, we use the following notation:
$p_{xy0|t} = P(X=x, Y=y, M=0\mid T=t), p_{+y1|t} = P(Y=y, M=1\mid T=t)$,
for $x=0, 1, \ldots, J-1$, $y=0, 1, \ldots, K-1$ and $t= 0, 1$.
All these probabilities are identifiable by the observed data.
\\

\noindent {\it Proof of Theorem 1.}
Because $P(M = m \mid T = t, Y = y)$ is identifiable,
$P(X = x, Y = y\mid T = t) = p_{xy0|t}  / P(M=0 \mid T = t, Y = y)$ is identifiable.
The joint distribution can be identified by
$p_{txym}= P(T = t) P(X = x, Y = y\mid T = t) P(M = m\mid T = t, Y = y)$.
\ \ \ $\Box$\\

\noindent {\it Proof of Theorem 2.}
(1)
Since $M\ind Y|(T,X)$, we have $P(Y = y\mid T = t, X = x) = P(Y = y\mid T = t, X = x, M=0)$.
Therefore, $P(Y = y \mid T = t, X = x)$ is identifiable, and so is $CE_x$ by the ignorability assumption.

\noindent (2) By $M\ind Y|(T,X)$, we obtain that $p_{xy0|t} = P(X = x, Y=y\mid T = t)  P(M=0\mid T = t, X = x)$ and thus
\begin{eqnarray*}
p_{+y1|t}    &=&  \sum\limits_{x=0}^{J-1}  P(X = x, Y = y, M=1\mid T = t) \\
                          &=& \sum\limits_{x=0}^{J-1}   P(X = x, Y=y\mid T= t) P(M=1\mid T = t, X = x)  \\
                          &=& \sum\limits_{x=0}^{J-1}   \frac{  P(M=1\mid T=t, X=x)  }{  P(M=0\mid T=t, X=x) } p_{xy0|t} \\
                      &=& \sum\limits_{x=0}^{J-1} \xi_{tx}   p_{xy0|t}  ,
\end{eqnarray*}
where $\xi_{tx} = P(M=1\mid T = t, X = x) / P(M=0\mid T = t, X = x )$.
It can be rewritten as
\begin{eqnarray}
\begin{pmatrix}
p_{000|t}& \cdots& p_{(J-1)00|t}\\
\vdots& &\vdots\\
p_{0(K-1)0|t}&\cdots &  p_{(J-1)(K-1)0|t}
\end{pmatrix}
\begin{pmatrix}
\xi_{t0}\\
\vdots\\
\xi_{t(J-1)}
\end{pmatrix}
=
\begin{pmatrix}
 p_{+01|t}\\
\vdots\\
 p_{+(K-1)1|t}
\end{pmatrix}. \label{eq::model1_1}
\end{eqnarray}
Notice that the probabilities in both sides of the equation are identifiable.
Thus the parameters $(\xi_{t0}, \ldots, \xi_{t(J-1)} )$ in (\ref{eq::model1_1}) are identifiable if
the solution of the linear equations (\ref{eq::model1_1}) is unique, which is equivalent to
Rank$(\bm{\Theta}_t ) = J$, where $\bm{\Theta}_t = (p_{txy0})$.
By the definition of $\xi_{tx} $, the missing data mechanism $P(M=0\mid T=t, X=x) = 1/(1+\xi_{tx})$ can be identified after identifying $\xi_{tx}$.
The conditional distribution of $(X,Y)$ given $T$ can be identified by
$
P(X = x, Y = y\mid T = t) =  p_{xy0|t} / P(M=0\mid T = t, X = x),
$
and the identifiability of $P(T=t)$ is obvious. The joint distribution can be identified by
$p_{txym} = P(T = t)P(X = x, Y = y\mid T = t) P(M = m\mid T = t, X = x)$.

\noindent (3) When $X$ is binary, Rank$(\bm{\Theta}_t ) = J=2$ holds if and only if there exists $y\neq 0$, such that $ p_{000|1}  p_{100|1}
 \not =
  p_{0y0|1} p_{1y0|1} $,
 \begin{eqnarray}
 \frac{P(X=0, Y=0, M=0\mid T=t)}{ P(X=0, Y=y, M=0\mid T=t) } \neq  \frac{P(X=1, Y=0, M=0\mid T=t)}{ P(X=1, Y=y, M=0\mid T=t) }  .\label{eq::rank2}
 \end{eqnarray}
By $P(X=x, Y=y, M=0\mid T=t) = P(X=x, Y=y\mid T=t) P(M=0\mid T=t, X=x)$, (\ref{eq::rank2}) is equivalent to
  \begin{eqnarray*}
 \frac{P(X=0, Y=0\mid T=t)}{ P(X=0, Y=y\mid T=t) } \neq  \frac{P(X=1, Y=0\mid T=t)}{ P(X=1, Y=y\mid T=t) }, \text{ i.e., }X\nind Y|(T=t).
 \end{eqnarray*}
By $P(X=x, Y=y, M=0\mid T=t) = P(M=0\mid T=t) P(X=x, Y=y\mid T=t, M=0)$, (\ref{eq::rank2}) is equivalent to
  \begin{eqnarray*}
 \frac{P(X=0, Y=0\mid T=t, M=0)}{ P(X=0, Y=y\mid T=t, M=0) } \neq  \frac{P(X=1, Y=0\mid T=t, M=0)}{ P(X=1, Y=y\mid T=t,M=0) } , \text{ i.e., } X\nind Y|(T=t, M=0).
 \end{eqnarray*}
\ \ \  $\Box$ \\

\noindent {\it Proof of Theorem 3.}
(1)
By $M \ind T|(X,Y)$, we have
\begin{eqnarray}\label{eq::thm3_1}
\frac{ p_{x10|1}   }{p_{x10|0}   }   = \frac{  P(Y= 1\mid T=1, X= x) P(X=x\mid T=1)}{P(Y = 1\mid T=0, X = x)P(X=x\mid T=0)}
\end{eqnarray}
and
\begin{eqnarray}\label{eq::thm3_2}
\frac{ p_{x00|1}   }{p_{x00|0}   }   = \frac{  P(Y = 0\mid T=1, X= x) P(X=x\mid T=1)}{  P(Y=0\mid T=0, X=x)P(X=x\mid T=0)}.
\end{eqnarray}
Dividing (\ref{eq::thm3_1}) by (\ref{eq::thm3_2}), we obtain that
$$
COR_x
=  \frac{  P(Y= 1\mid T=1, X= x)  P(Y=0\mid T=0, X=x) }{P(Y = 1\mid T=0, X = x)   P(Y = 0\mid T=1, X=x)}
=  \frac{ p_{x10|1}  p_{x00|0}  }{p_{x10|0}   p_{x00|1} }
$$
are identifiable.
Since $COR_x >1$, $=1$ and $<1$ are equivalent to
$CE_x >0$, $=0$ and $<0$,
the signs of $CE_x$ are identifiable.

\noindent (2)
By $M \ind T|(X,Y)$ and $T\ind X$, we have
$$
\frac{ p_{x10|1}   }{p_{x10|0}   }   = \frac{  P(Y= 1\mid T=1, X= x)}{P(Y = 1\mid T=0, X = x)} \mbox{ and }
\frac{ p_{x00|1}   }{p_{x00|0}   }   = \frac{  1 - P(Y = 1\mid T=1, x)}{ 1 -  P(Y=1\mid T=0, X=x)}.
$$
Therefore, we obtained the following linear equations for $P(Y=1\mid T=1, X=x)  $ and $P(Y=1\mid T=0, X=x)$:
\begin{eqnarray}
\begin{pmatrix}
p_{x10|0}   & p_{x10|1}  \\
p_{x00|0}   &  p_{x00|1}
\end{pmatrix}
\begin{pmatrix}
P(Y=1\mid T=1, X=x)\\
 -P(Y=1\mid T=0, X=x)
\end{pmatrix}
=
\begin{pmatrix}
0\\
p_{x00|0}  - p_{x00|1}
\end{pmatrix}.\label{eq::3}
\end{eqnarray}
When $ p_{x10|0}  p_{x00|1}    \not =  p_{x10|1}   p_{x00|0} $, the solution of (\ref{eq::3}) is unique, and thus $P(Y=1\mid T=1, X=x)  $ and $P(Y=1\mid T=0, X=x)  $ are identifiable. When $ p_{x10|0}  p_{x00|1}    =  p_{x10|1} p_{x00|0} $, we have $Y\ind T|X$ from $M \ind T|(X,Y)$, and thus we get $P(Y=1\mid T=1, X=x)  = P(Y=1\mid T=0, X=x)$ and $CE_x= 0$. So the causal effects $CE_x$ are identifiable.

\noindent (3)
By the definition of $p_{xy0|t}  $ and $p_{+y1|t}  $, we obtain that $p_{xy0|t}  = P(X=x, Y=y\mid T=t) P(M=0\mid X=x, Y=y)$  and
that for $J=2$,
\begin{eqnarray*}
p_{+y1|t}    &=&  \sum\limits_{x=0,1}  P(X=x,  Y=y, M=1\mid T=t) \\
                          &=& \sum\limits_{x=0,1}  P(X=x, Y=y \mid T=t) P(M=1 \mid X=x, Y=y)  \\
                          &=& \sum\limits_{x=0, 1}   \frac{  P(M=1 \mid X=x, Y=y) }{  P(M=0 \mid X=x, Y=y) } p_{xy0|t}\\
                          &=&  \sum\limits_{x=0,1}  \kappa_{xy} p_{xy0|t},
\end{eqnarray*}
where $ \kappa_{xy} =  P(M=1 \mid X=x, Y=y)  / P(M=0 \mid X=x, Y=y).$
Thus we have
\begin{eqnarray}
\begin{pmatrix}
p_{0y0|1}  &  p_{1y0|1} \\
p_{0y0|0}  &  p_{1y0|0}
\end{pmatrix}
\begin{pmatrix}
\kappa_{0y}\\
\kappa_{1y}
\end{pmatrix}
=
\begin{pmatrix}
p_{+y1|1} \\
p_{+y1|0} \\
\end{pmatrix}, y=0,1, \ldots, K-1. \label{eq::model3_1_X01}
\end{eqnarray}
The parameters $(\kappa_{0y}, \kappa_{1y})$ are identifiable
if and only if  $ p_{0y0|1}  p_{1y0|0}   \not=   p_{1y0|1} p_{0y0|0}  $.
After identifying $\kappa_{xy} $, the missing data mechanism can be identified by $P(M=0\mid X=x, Y=y) = 1/(1 + \kappa_{xy} )$,
and then
$P(X=x, Y=y\mid T=t) =  p_{xy0|t}   / P(M=0\mid X=x, Y=y)$ is identifiable.
Finally we identify $p_{txym} = P(T=t) P(X=x, Y=y\mid T=t) P(M=m\mid X=x,Y=y) $.

Further, the condition $ p_{0y0|1}  p_{1y0|0}   \not=   p_{1y0|1} p_{0y0|0}  $ is
\begin{eqnarray}
\frac{P(X=0, Y=y, M=0\mid T=1)}{P(X=0, Y=y, M=0\mid T=0)} \neq \frac{P(X=1, Y=y, M=0\mid T=1)}{P(X=1, Y=y, M=0\mid T=0)} .\label{eq::rank3}
\end{eqnarray}
By $P(X=x, Y=y, M=0\mid T=t) = P(Y=y\mid T=t)P(X=x\mid T=t, Y=y)P(M=0\mid X=x, Y=y) $, (\ref{eq::rank3}) is equivalent to
\begin{eqnarray*}
\frac{P(X=0\mid T=1, Y=y)}{P(X=0\mid T=0, Y=y)} \neq \frac{P(X=1\mid T=1, Y=y)}{P(X=1\mid T=0, Y=y)}, \text{ i.e., } X\nind T|(Y=y).
\end{eqnarray*}
By $P(X=x, Y=y, M=0\mid T=t) = P(Y=y, M=0\mid T=t)P(X=x\mid T=t, Y=y, M=0) $, (\ref{eq::rank3}) is equivalent to
\begin{eqnarray*}
\frac{P(X=0\mid T=1, Y=y, M=0)}{P(X=0\mid T=0, Y=y, M=0)} \neq \frac{P(X=1\mid T=1, Y=y, M=0)}{P(X=1\mid T=0, Y=y, M=0)}, \text{ i.e., } X\nind T|(Y=y, M=0).
\end{eqnarray*}
\ \ \ $\Box$
\\

\noindent {\it Proof of Corollary 1.}
The first part of the corollary is obvious from Theorem 2.
By $M\ind (T,Y)|X$, we have $\xi_{tx} = \kappa_{xy} =   P(M=1\mid X=x) / P(M=0 \mid X=x)$, denoted as $\gamma_x$.
From (\ref{eq::model1_1}) for $J=2$,
the solution of $(\gamma_1, \gamma_0)$ is unique if there exists a $t\in \{0,1\}$ such that $X\nind Y | (T=t)$. From (\ref{eq::model3_1_X01}), the solution of $(\gamma_1, \gamma_0)$ is unique if there exists a $y\in \{0,...,K-1\}$ such that $X\nind T|(Y=y)$.
Since $X\ind (T,Y)$ is equivalent to $X\ind T|Y$ and $X\ind Y|T$, we have that $X\nind (T,Y)$ is equivalent to $X\nind T|Y$ or $X\nind Y|T$.
Therefore, $(\gamma_1, \gamma_0)$ is identifiable if $X\nind (T,Y)$, which is equivalent to $X\nind (T,Y)|(M=0)$, 
since all the conditions can be replaced by further conditioning on $M=0$ as proved in Theorem 2 and 3.
After identifying $(\gamma_1, \gamma_0)$,
we can show the identifiability of the joint distribution
in the same way as the proof of Theorem 2.
\ \ \ $\Box$
\\

\noindent {\it Proof of Theorem 4.}
By $p_{xy0|t} = P(X = x\mid T = t) P(Y = y\mid T = t, X = x) P(M=0\mid T = t, X = x, Y = y)$, we have
\begin{eqnarray*}
p_{+y1|t}    &=&  \sum\limits_{x=0,1} P(X = x, Y = y, M=1\mid T = t) \\
                          &=& \sum\limits_{x=0,1} P(X = x\mid T = t) P(Y = y\mid T = t, X = x) P(M=1 \mid T = t, X = x, Y = y)  \\
                          &=& \sum\limits_{x=0,1} \frac{  P(M=1 \mid T = t, X = x, Y = y) }{  P(M=0 \mid T = t, X = x, Y = y) } p_{xy0|t}.
\end{eqnarray*}
By the logistic missing data mechanism, we get
\begin{equation}\label{eq::moment4}
p_{+y1|t} = \sum\limits_{x=0,1}  \exp\{\beta_0+ \beta_T t + \beta_X x + \beta_Y y  \} p_{xy0|t}.
\end{equation}
Equation (\ref{eq::moment4}) implies the following equations
 \begin{numcases}{}
p_{+01|0} = p_{000|0}  \exp(\beta_0) +   p_{100|0}  \exp(\beta_0 + \beta_X),   \label{eq::model4::1}\\
p_{+01|1}  = p_{000|1}  \exp(\beta_0 + \beta_T )  +  p_{100|1}  \exp(\beta_0 + \beta_T +\beta_X), \label{eq::model4::2} \\
p_{+11|0}  = p_{010|0}   \exp(\beta_0 + \beta_Y) + p_{110|0}  \exp(\beta_0 + \beta_X + \beta_Y), \label{eq::model4::3}\\
p_{+11|1} = p_{010|1}  \exp(\beta_0 + \beta_T + \beta_Y)  + p_{110|1}  \exp(\beta_0 + \beta_T +\beta_X + \beta_Y).  \label{eq::model4::4}
 \end{numcases}
Let $ A = \exp(\beta_0)$, $B =  \exp(\beta_X)$, $C =    \exp(\beta_0 + \beta_T )$ and $D =   \exp(\beta_0 + \beta_Y)$, and thus $CD/A =  \exp(\beta_0 + \beta_T + \beta_Y).$
From (\ref{eq::model4::1}) to (\ref{eq::model4::4}),
we have $A = p_{+01|0} / (p_{000|0}  +  p_{100|0} B)$,
$C = p_{+01|1} / (p_{000|1} + p_{100|1}  B)$,
$D = p_{+11|0} / (p_{010|0} + p_{110|0}  B)$ and
$CD/A = p_{+11|1} / (p_{010|1}   + p_{110|1} B)$ respectively.
By $CD = A(CD/A)$, we get the following quadratic equation of $B$
\begin{eqnarray*}
\frac{   p_{+01|1}  p_{+11|0}    }{(  p_{000|1}    +   p_{100|1}  B)(p_{010|0}    + p_{110|0}  B) }
=
\frac{   p_{+01|0} p_{+11|1}    }{( p_{000|0}  +  p_{100|0}  B )(p_{010|1}   + p_{110|1} B)},
\end{eqnarray*}
which is equivalent to
$EB^2 + FB + G = 0$,
where
\begin{eqnarray*}
E &=& \frac{   p_{110|0} p_{100|1} }{    p_{+11|0}  p_{+01|1}  }
-
\frac{  p_{100|0} p_{110|1}  }{   p_{+01|0}  p_{+11|1}    },  \label{eq::E}\\
F &=& \frac{p_{110|0} p_{000|1}  + p_{010|0} p_{100|1}  }
{   p_{+11|0}  p_{+01|1}    }
-
\frac{  p_{100|0} p_{010|1}  + p_{000|0} p_{110|1} }
{     p_{+01|0} p_{+11|1}   },\\
G &=& \frac{  p_{010|0} p_{000|1}  }{p_{+11|0}  p_{+01|1}   }
-
\frac{   p_{000|0} p_{010|1}  }{  p_{+01|0}  { p_{+11|1}} }. \label{eq::G}
\end{eqnarray*}
It is known that $B$ has only one positive solution if and only if $EG \leq 0$.
By algebraic operations, we know that $EG$ and $\{ OR_{YT|(M=1)}  - OR_{YT|(X=1,M=0)} \}
\{ OR_{YT|(M=1)}  - OR_{YT|(X=0,M=0)} \}$ have the same sign.
Thus $OR_{YT|(M=1)} $ is between $ OR_{YT|(X=1,M=0)} $
and $ OR_{YT|(X=0,M=0)} $ if and only if $EG\leq0$
If $B$ has only a positive solution,
$\beta_{0}$, $\beta_T$, $\beta_{X}$ and $\beta_Y$ are identifiable.
After identifying the missing data mechanism, we can identify $P(X=x,Y=y\mid T=t) = p_{xy0|t}  \left[   1 + \exp\{\beta_0 + \beta_T t + \beta_{X} x +  \beta_{Y}t \}   \right]$, and thus
the joint distribution can be identified by $p_{txym}=P(T = t) P(X = x,Y = y\mid T = t) P(M = m\mid T = t, X = x, Y = y)$.
\ \ \ $\Box$ \\

\noindent {\it Proof of Theorem 5.}
In $P(Y=1 \mid T=1, X = x) = (p_{1x10} + p_{1x11}) / (p_{1x00} + p_{1x10} + p_{1x01} + p_{1x11})$,
only $p_{1x01}$ and $p_{1x11}$ cannot be identified,
but they are subject to the constraints
$0\leq p_{1x01}\leq p_{1+01}$ and $0\leq p_{1x11} \leq p_{1+11}$.
It can be seen from the equation that $P(Y=1\mid T=1, X=x)$ is an decreasing function of $p_{1x01}$ and an increasing function of $ p_{1x11}$.
Therefore, the minimum of $P(Y=1\mid T=1, X=x)$ is $p_{1x10} / (p_{1x00} + p_{1x10} + p_{1+01})$, which is obtained at $p_{1x01} = p_{1+01}$ and $p_{1x11} = 0$; and the maximum of $P(Y=1\mid T=1, X=x)$ is $(p_{1x10} + p_{1+11}) / (p_{1x00} + p_{1x10} + p_{1+11})$, which is obtained at
$p_{1x01} = p_{1+01}$ and $p_{1x11} = 0$.

Similarly, $P(Y=1\mid T=0, X=x) = (p_{0x10} + p_{0x11}) / (p_{0x00} + p_{0x10} + p_{0x01} + p_{0x11})$, which is an decreasing function of $p_{0x01}$ and an increasing function of $ p_{0x11}$.
The minimum of $ P(Y=1\mid T=0, X=x)$ is $(p_{0x10}) / (p_{0x00} + p_{0x10}   + p_{0+01})$ at $p_{0x01} = p_{0+01}$ and $p_{0111} = 0$, and the maximum of $P(Y=1\mid T=0, X=x)$ is $(p_{0x10} + p_{0+11}) / (p_{0x00} + p_{0x10} + p_{0+11})$ at $p_{0x01} = 0$ and $p_{0x11 } = p_{0+11}$.
Since $\mathcal{D}[p_1, p_0]$ satisfies $\partial\mathcal{D}/\partial p_1 >0$ and $\partial\mathcal{D}/\partial p_0 <0$, we proved this theorem.
\ \ \ $\Box$

\section{Table of the Simulation Studies}

We generated observed data sets from all missing mechanisms,
and applied five models based on five missing mechanisms to each data set,
denoted as ``$M_1$'' to ``$M_5$'' in Table \ref{tb::simulation}.
Thus we also show the sensitivities
when the missing mechanism is not correctly specified.
Applying various models to each data set,
we estimated $CE_0 = \log(COR_0)$ and $CE_1 = \log(COR_1)$ using the MLEs and the posterior medians,
and we calculated the bounds of mechanism 5 given in Theorem 5.
Using the Gibbs samplers to get the Bayesian credible intervals is more direct than using the likelihood-based methods.
Gibbs samplers were run $10000$ times with burning in after the $5000$-th iteration.
We did the simulation studies under sample sizes $500$ and $1000$, and the processes were repeated $1000$ times.
Table \ref{tb::simulation} presents the results of the simulation studies, which contain the average biases (bias$_{\text{EM}}$ and bias$_{\text{Gibbs}}$) and the mean square errors (MSE$_{\text{EM}}$ and MSE$_{\text{Gibbs}}$),
the coverage proportions of the $95\%$ posterior credible intervals from the Gibbs sampler (CP$_{\text{Gibbs}}$), and the means of the upper and lower bounds obtained from mechanism 5 (upper and lower).
The data generating processes $M_1$, $M_2$, $M_3$ and $M_4$ satisfy the specification conditions of the missing mechanisms 1 to 4, respectively. Therefore, the first four diagonal blocks in boldface of the upper and the lower panels of Table \ref{tb::simulation} show very good performances of our proposed methods with small average biases, small mean square errors and reasonable coverage proportions, if the model specifications are correct.
And the mean square errors decrease as the sample sizes increase under correctly specified models.
However, misspecification of the missing data mechanisms can cause severe biases and poor coverage proportions, which suggests the importance of the specification of the missing data mechanism.
For the process `` $M^*$'', we first selected a missing mechanism from the mechanisms 1 to 4 based on the log likelihood functions because all these mechanisms have the same numbers of parameters, and then we apply the selected mechanism to the data set.
The rows labeled ``$M^*$'' in Table \ref{tb::simulation} are the average biases and mean square errors of the estimates based on the process,
where there are a fewer number of large biases and MSEs than
the other models
and the largest bias and MSE are less than those of the other models.
Although mechanism 5 does not need any assumption of missing mechanism,
the bounds are too wide and all of them cover zero.

\newsavebox{\tablebox}
\begin{lrbox}{\tablebox}
\begin{tabular}{rrrrrrrrrrrrrrrr}
\hline
 &  & \multicolumn{2}{c}{$M_1$} &  & \multicolumn{2}{c}{$M_2$}&  & \multicolumn{2}{c}{$M_3$} & & \multicolumn{2}{c}{$M_4$}& & \multicolumn{2}{c}{$M_5$} \\
 \cline{3-4} \cline{6-7} \cline{9-10} \cline{12-13} \cline{15-16}
 & &$CE_0$ &  $CE_1$ & &$CE_0$ &$CE_1$ &&$CE_0$ &$CE_1$ &&$CE_0$ &$CE_1$&&$CE_0$ &$CE_1$ \\
 \hline
$N=500$\\
$M_1$ & bias$_{\text{EM}}$ &{\bfseries 0.046 }& {\bfseries -0.040} &  & 0.343 & 0.239 &  & -0.596 & -0.714 &  & 0.253 & 0.129 &  & 0.144 & 0.216 \\
 & MSE$_{\text{EM}}$ &{\bfseries  0.152} & {\bfseries 0.121} &  & 0.395 & 0.153 &  & 0.469 & 0.739 &  & 0.265 & 0.126 &  & 0.230 & 0.150 \\
 & bias$_{\text{Gibbs}}$ &{\bfseries  0.038} & {\bfseries -0.041} &  & 0.332 & 0.245 &  & -0.598 & -0.694 &  & 0.243 & 0.131 &  & 0.136 & 0.219 \\
 & MSE$_{\text{Gibbs}}$ & {\bfseries 0.149} & {\bfseries 0.121} &  & 0.381 & 0.156 &  & 0.471 & 0.704 &  & 0.256 & 0.127 &  & 0.226 & 0.153 \\
 & CP$_{\text{Gibbs}}$ & {\bfseries 0.956} & {\bfseries 0.949 } &  & 0.889 & 0.869 &  & 0.566 & 0.702 &  & 0.924 & 0.924 &  & 0.944 & 0.888 \\
$M_2$ & bias$_{\text{EM}}$ & -0.506 & -0.718 &  &{\bfseries  0.073} &{\bfseries  -0.025 }&  & -0.161 & -0.593 &  & 0.133 & -0.193 &  & 0.180 & 0.124 \\
 & MSE$_{\text{EM}}$ & 0.438 & 0.688 &  & {\bfseries 0.290} & {\bfseries 0.127} &  & 0.256 & 0.597 &  & 0.334 & 0.219 &  & 0.276 & 0.206 \\
 & bias$_{\text{Gibbs}}$ & -0.537 & -0.691 &  & {\bfseries -0.010 }& {\bfseries 0.004} &  & -0.220 & -0.486 &  & 0.059 & -0.150 &  & 0.113 & 0.170 \\
& MSE$_{\text{Gibbs}}$ & 0.428 & 0.648 &  & {\bfseries 0.290} & {\bfseries 0.131} &  & 0.302 & 0.489 &  & 0.307 & 0.201 &  & 0.268 & 0.211 \\
 & CP$_{\text{Gibbs}}$ & 0.684 & 0.593 &  & {\bfseries 0.935} & {\bfseries 0.947 }&  & 0.867 & 0.844 &  & 0.912 & 0.918 &  & 0.942 & 0.898 \\
$M_3$ & bias$_{\text{EM}}$ & -0.130 & -0.842 &  & 0.237 & -0.362 &  & {\bfseries 0.037 }& {\bfseries -0.053 }&  & 0.248 & -0.566 &  & 0.273 & -0.539 \\
& MSE$_{\text{EM}}$ & 0.131 & 0.869 &  & 0.290 & 0.234 &  & {\bfseries 0.153} & {\bfseries 0.222 }&  & 0.205 & 0.464 &  & 0.235 & 0.435 \\
& bias$_{\text{Gibbs}}$ & -0.125 & -0.796 &  & 0.255 & -0.320 &  & {\bfseries -0.033 }& {\bfseries 0.048 }&  & 0.268 & -0.545 &  & 0.302 & -0.526 \\
& MSE$_{\text{Gibbs}}$ & 0.133 & 0.794 &  & 0.288 & 0.209 &  & {\bfseries 0.159 }& {\bfseries 0.196 }&  & 0.216 & 0.444 &  & 0.241 & 0.414 \\
& CP$_{\text{Gibbs}}$ & 0.925 & 0.454 &  & 0.920 & 0.844 &  &{\bfseries  0.942 }&{\bfseries  0.965 }&  & 0.913 & 0.673 &  & 0.908 & 0.694 \\
$M_4$ & bias$_{\text{EM}}$ & -0.395 & -0.594 &  & 0.154 & 0.003 &  & -0.029 & -0.026 &  & {\bfseries 0.027 }& {\bfseries -0.047} &  & -0.073 & 0.093 \\
& MSE$_{\text{EM}}$ & 0.298 & 0.527 &  & 0.292 & 0.131 &  & 0.152 & 0.298 &  & {\bfseries 0.236} & {\bfseries 0.143 }&  & 0.256 & 0.135 \\
& bias$_{\text{Gibbs}}$ & -0.505 & -0.599 &  & 0.087 & 0.035 &  & -0.040 & 0.126 &  &{\bfseries  -0.048} & {\bfseries -0.033} &  & -0.146 & 0.114 \\
& MSE$_{\text{Gibbs}}$ & 0.417 & 0.535 &  & 0.283 & 0.143 &  & 0.156 & 0.307 &  & {\bfseries 0.251} & {\bfseries 0.143} &  & 0.277 & 0.142 \\
& CP$_{\text{Gibbs}}$ & 0.679 & 0.689 &  & 0.942 & 0.938 &  & 0.941 & 0.786 &  & {\bfseries 0.935} & {\bfseries 0.945} &  & 0.916 & 0.938 \\
 $M^*$ & bias$_{\text{EM}}$ & -0.075 & -0.158 &  & 0.207 & 0.103 &  & -0.222 & -0.300 &  & 0.098 & 0.027 &  & -0.031 & 0.128 \\
& MSE$_{\text{EM}}$ & 0.180 & 0.194 &  & 0.331 & 0.138 &  & 0.259 & 0.486 &  & 0.246 & 0.125 &  & 0.258 & 0.133 \\
$M_5$ & upper & 4.094 & 1.459 &  & 5.204 & 2.040 &  & 4.069 & 2.465 &  & 4.658 & 1.667 &  & 4.900 & 1.993 \\
& lower & -0.225 & -2.690 &  & -0.770 & -2.539 &  & -0.052 & -3.407 &  & -0.603 & -2.525 &  & -0.995 & -2.651 \\
 $N=1000$\\
$M_1$ & bias$_{\text{EM}}$ & {\bfseries 0.028 }& {\bfseries -0.009} &  & 0.290 & 0.243 &  & -0.611 & -0.683 &  & 0.275 & 0.149 &  & 0.123 & 0.221 \\
& MSE$_{\text{EM}}$ & {\bfseries 0.070} &{\bfseries  0.052 }&  & 0.199 & 0.103 &  & 0.430 & 0.579 &  & 0.173 & 0.073 &  & 0.114 & 0.099 \\
 & bias$_{\text{Gibbs}}$ & {\bfseries 0.024 }& {\bfseries -0.009 }&  & 0.286 & 0.246 &  & -0.611 & -0.674 &  & 0.270 & 0.150 &  & 0.120 & 0.222 \\
 & MSE$_{\text{Gibbs}}$ & {\bfseries 0.070 }& {\bfseries 0.052 }&  & 0.196 & 0.104 &  & 0.431 & 0.566 &  & 0.170 & 0.074 &  & 0.113 & 0.100 \\
& CP$_{\text{Gibbs}}$ & {\bfseries 0.954 }& {\bfseries 0.966 }&  & 0.877 & 0.807 &  & 0.294 & 0.448 &  & 0.865 & 0.892 &  & 0.938 & 0.847 \\
$M_2$ & bias$_{\text{EM}}$ & -0.552 & -0.680 &  &{\bfseries  0.034} &{\bfseries  -0.012 }&  & -0.186 & -0.598 &  & 0.206 & -0.155 &  & 0.174 & 0.133 \\
& MSE$_{\text{EM}}$ & 0.397 & 0.536 &  & {\bfseries 0.121 }&{\bfseries  0.058} &  & 0.148 & 0.467 &  & 0.171 & 0.113 &  & 0.144 & 0.106 \\
& bias$_{\text{Gibbs}}$ & -0.581 & -0.667 &  & {\bfseries -0.006 }& {\bfseries 0.001 }&  & -0.212 & -0.551 &  & 0.162 & -0.129 &  & 0.143 & 0.168 \\
& MSE$_{\text{Gibbs}}$ & 0.410 & 0.518 &  &{\bfseries  0.124 }& {\bfseries 0.059 }&  & 0.166 & 0.419 &  & 0.167 & 0.108 &  & 0.140 & 0.122 \\
 & CP$_{\text{Gibbs}}$ & 0.456 & 0.306 &  & {\bfseries 0.954 }& {\bfseries 0.954 }&  & 0.848 & 0.603 &  & 0.902 & 0.898 &  & 0.937 & 0.884 \\
$M_3$ & bias$_{\text{EM}}$ & -0.140 & -0.802 &  & 0.172 & -0.352 &  & {\bfseries 0.017} & {\bfseries -0.020 }&  & 0.266 & -0.531 &  & 0.250 & -0.536 \\
 & MSE$_{\text{EM}}$ & 0.074 & 0.710 &  & 0.116 & 0.168 &  & {\bfseries 0.073 }& {\bfseries 0.106 }&  & 0.142 & 0.350 &  & 0.134 & 0.350 \\
& bias$_{\text{Gibbs}}$ & -0.136 & -0.779 &  & 0.188 & -0.334 &  & {\bfseries -0.016} & {\bfseries 0.060 }&  & 0.278 & -0.525 &  & 0.263 & -0.533 \\
 & MSE$_{\text{Gibbs}}$ & 0.074 & 0.674 &  & 0.122 & 0.157 &  & {\bfseries 0.074} & {\bfseries 0.115 }&  & 0.148 & 0.345 &  & 0.141 & 0.348 \\
& CP$_{\text{Gibbs}}$ & 0.915 & 0.122 &  & 0.922 & 0.702 &  & {\bfseries 0.949 }& {\bfseries 0.950 }&  & 0.849 & 0.439 &  & 0.857 & 0.427 \\
$M_4$ & bias$_{\text{EM}}$ & -0.420 & -0.557 &  & 0.106 & 0.017 &  & -0.060 & -0.010 &  & {\bfseries 0.058 }& {\bfseries -0.014 }&  & -0.079 & 0.101 \\
& MSE$_{\text{EM}}$ & 0.240 & 0.387 &  & 0.125 & 0.061 &  & 0.075 & 0.154 &  & {\bfseries 0.112} & {\bfseries 0.063 }&  & 0.121 & 0.068 \\
& bias$_{\text{Gibbs}}$ & -0.475 & -0.560 &  & 0.079 & 0.032 &  & -0.057 & 0.127 &  & {\bfseries 0.024} & {\bfseries -0.007 }&  & -0.114 & 0.111 \\
& MSE$_{\text{Gibbs}}$ & 0.296 & 0.391 &  & 0.122 & 0.064 &  & 0.075 & 0.189 &  & {\bfseries 0.112 }& {\bfseries 0.064} &  & 0.130 & 0.071 \\
& CP$_{\text{Gibbs}}$ & 0.608 & 0.492 &  & 0.948 & 0.949 &  & 0.948 & 0.810 &  & {\bfseries 0.944 }& {\bfseries 0.958 }&  & 0.933 & 0.931 \\
$M^*$ & bias$_{\text{EM}}$ & -0.044 & -0.079 &  & 0.169 & 0.121 &  & -0.214 & -0.249 &  & 0.115 & 0.040 &  & -0.051 & 0.124 \\
& MSE$_{\text{EM}}$ & 0.087 & 0.093 &  & 0.160 & 0.078 &  & 0.187 & 0.301 &  & 0.127 & 0.058 &  & 0.122 & 0.072 \\
$M_5$ & upper & 4.078 & 1.473 &  & 5.133 & 2.024 &  & 4.044 & 2.457 &  & 4.683 & 1.677 &  & 4.862 & 1.989 \\
 & lower & -0.223 & -2.649 &  & -0.764 & -2.522 &  & -0.048 & -3.371 &  & -0.584 & -2.492 &  & -0.970 & -2.642 \\ \hline
\end{tabular}
\end{lrbox}

\begin{table}[htb]
\begin{center}
\caption{Simulation Studies. The true values are $CE_0= 2.773$ and $CE_1 = -0.847$.
The upper and lower panels shows the results of sample sizes $500$ and $1000$.
The columns with ``$M_1$'' to ``$M_5$'' correspond to the five data generating processes in Section 6. The rows with ``$M_1$'' to ``$M5$'' are the results obtained by the methods under mechanisms 1 to 5, and the rows with ``$M^*$'' are the results obtained by the selected mechanism based on the likelihoods.}\label{tb::simulation}
\resizebox{\textwidth}{!}{\usebox{\tablebox}}
\end{center}
\end{table}

\section{Balance Checking for Covariates in the Job Training Data}

The results are shown in Table \ref{tb::balance}.

\begin{table}[ht]
\begin{center}
\caption{Balancing Checking of the Covariates}\label{tb::balance}
\begin{tabular}{rlrrr}
  \hline
 & variables & t statistic & $95\%$ confidence interval & $p$-value \\ 
  \hline
1 & age & 1.114 &  $[ -0.583 , 2.108]$ & 0.266 \\ 
  2 & education & 1.442 & $[-0.094 , 0.609]$ & 0.150 \\ 
  3 & black & 0.458 & $[ -0.054 , 0.086]$ & 0.647 \\ 
  4 & hispanic & -1.857 & $[ -0.099 , 0.003]$ & 0.064 \\ 
  5 & no degree & -3.108 & $[-0.207 , -0.046]$ & 0.002$^*$ \\ 
  6 & married & 0.967 & $[-0.037 , 0.107]$ & 0.334 \\ 
  7 & no job in 1974 & -0.975 & $[-0.126 , 0.043]$ & 0.330 \\ 
  8 & no job in 1975 & -1.830 &$[ -0.176 , 0.006]$ & 0.068 \\ 
   \hline
\end{tabular}
\end{center}
\end{table}

\end{document}